\documentclass[10pt, reqno]{amsart}
\usepackage{amssymb,amsmath,amscd,amsfonts,amsthm,mathrsfs}
\usepackage[french, english]{babel}
\usepackage{verbatim}
\usepackage{graphicx}             
\input xy
\xyoption{all}
\usepackage[all]{xy}
\usepackage{color}

\usepackage[latin1]{inputenc}

\newcommand{\N}{\mathbb{N}}

\newcommand{\C}{\mathbb{C}}
\newcommand{\Ac}{\mathcal{A}}
\newcommand{\Cc}{\mathcal{C}}
\newcommand{\Ec}{\mathcal{E}}
\newcommand{\Fc}{\mathcal{F}}
\newcommand{\Ic}{\mathcal{I}}
\newcommand{\Kc}{\mathcal{K}}
\newcommand{\Hc}{\mathcal{H}}
\newcommand{\Sc}{\mathcal{S}}
\newcommand{\Tc}{\mathcal{T}}
\newcommand{\Vc}{\mathcal{V}}

\def\build#1_#2^#3{\mathrel{\mathop{\kern 0pt#1}\limits_{#2}^{#3}}}

\newtheorem{theorem}{Theorem}[section]
\newtheorem{propr}[theorem]{Proposition}
\newtheorem{lem}[theorem]{Lemma}
\newtheorem{rem}[theorem]{Remark}
\newtheorem{exa}[theorem]{Example}
\newtheorem{defi}[theorem]{Definition}
\newtheorem{coro}[theorem]{Corollary}
\newtheorem{thm}[theorem]{Theorem}
\newenvironment{dem}{\noindent {\it Proof: }
    \begin{quotation}\noindent}{\end{quotation}\hfill$\square$}

\numberwithin{equation}{section}

\title[Laplacian of a 2-simplicial complex]{The Discrete Laplacian of a 2-simplicial complex}

\begin{document}
\title[Laplacian of a 2-simplicial complex]{The Discrete Laplacian of a 2-simplicial complex}

\author{Yassin CHEBBI}
\address{Laboratoire de Math\'ematiques Jean Leray, Facult\'{e} des Sciences, CNRS,
 Universit\'{e} de Nantes, BP 92208, 44322 Nantes, France}
\address{Unit\'e de recherches Math\'ematiques et Applications, UR/13 ES 47,
 Facult\'e des Sciences de Bizerte de l'Universit\'e de Carthage, 7021-Bizerte, Tunisie}
\email{chebbiyassin88@gmail.com}
\email{yassin.chebbi@univ-nantes.fr}
\date{Version of \today}

\begin{abstract}
In this paper, we introduce the notion of oriented faces especially triangles
in a connected oriented locally finite graph. This framework then permits
to define the Laplace operator on this structure of the 2-simplicial complex.
We develop the notion of $\chi$-completeness for the graphs, based on the cut-off functions.
Moreover, we study essential self-adjointness of the discrete Laplacian
from the $\chi$-completeness geometric hypothesis .

 \hspace{-0.5cm}
 {\sc R\'esum\'e.}
Dans cet article, nous introduisons la notion de faces orient\'{e}es et plus
particuli\`{e}ment de triangles dans un graphe connexe orient\'{e} localement fini.
Ce cadre permet alors de d\'{e}finir l'op\'{e}rateur de Laplace sur cette structure
d'un 2-complexe simplicial. Nous d\'{e}veloppons la notion de $\chi$-compl\'{e}tude pour les graphes,
bas\'{e}e sur les fonctions de coupure. De plus, nous \'{e}tudions le caract\`{e}re essentiellement
auto-adjoint du Laplacien discret \`{a} partir de l'hypoth\`{e}se g\'{e}om\'{e}trique $\chi$-compl\'{e}tude.
\end{abstract}

\maketitle

\tableofcontents

\maketitle

\section{Introduction}
The impact of the geometry on the essential self-adjointness of the Laplacians
is studied in many areas of mathematics on Riemannian manifolds; see (\cite{C},
\cite{EL},\cite{G},\cite{M1}) and also on one-dimensional simplicial complexes;
see (\cite{AT},\cite{CTT},\cite{FLW},\cite{HKMW},\cite{M2},\cite{T}). Indeed, Laplacians
on Riemannian manifolds and simplicial complexes share a lot of common elements. Despite of
this, various geometric notions such as distance and completeness in the Riemannian
framework have no immediate analogue in the discrete setting. Combinatorial Laplacians
were originally studied on graphs, beginning with \emph{Kirchhoff} and his study of
electrical networks \cite{K}. Simplicial complexes can be viewed as generalizations
of graphs, as from any graph, we can form its clique complex, a 2-simplicial complex
whose faces correspond to the cliques of the graph. In this article, we take a
connected oriented locally finite graph and we introduce the oriented faces especially
triangles in such a way that every face is a triangle, so we can regard it as a
two-dimensional simplicial complex. This work presents a more general framework for the Laplacians
defined in terms of the combinatorial structure of a simplicial complex. The main result
of this work gives a geometric hypothesis to ensure essential self-adjointness
for the discrete Laplacian. We develop the $\chi$-completeness hypothesis for triangulations.
This hypothesis on locally finite graphs covers many situations that have been already studied
in \cite{AT}. The authors prove that the $\chi$-completeness is satisfied by graphs which are
complete for some intrinsic metric, as defined in \cite{FLW} and \cite{HKMW}.

The paper is structured as follows: In the second section, we will first present
the basic concepts about graphs or rather one-dimensional simplicial complexes.
Next, we introduce the notion of oriented faces more particularly triangles where
all the faces are triangles. This special structure of 2-simplicial complex
is called triangulation. Without loss of generality, we can assume that every
triangle is a face for simplicity sake. So this permits to define the Gau\ss -Bonnet
operator $T=d+\delta$ acting on triplets of functions, 1-forms and 2-forms. After that,
we define the discrete Laplacian by $L:=T^2$ which admits a decomposition according to
the degree
$$L:=L_0\oplus L_1\oplus L_2.
$$

In the third and fourth sections, we study the closability of the operators
which are used in the following sections. Next, we get started with refer to
\cite{AT} for the notion $\chi$-completness of the graphs and we develop this
geometric hypothesis for the triangulations in Definition \ref{defi4.2}.
Moreover, we have developed it through optimal example of the "\emph{tiangular tree}"
to produce a concrete way to prove a triangulation which is not $\chi$-complete,
based on the offspring function, we refer here to \cite{BGJ} for this notion.

In the fifth section, we address the main results concerning essential self-adjointness
for $T$ and $L.$ In the case of complete manifolds, there is a result of \emph{Chernoff};
see \cite{C}, and we also have for the discrete setting; see \cite{AT}, which conclude
that the Dirac operator is essentially self-adjoint. As a result, they prove essential
self-adjointness of the Laplace-Beltrami operator. So, we take this idea to make the relationship
between $T$ and $L$ about the essential self-adjointness, when the triangulation is $\chi$-complete.

In the final section, we present a particular example of a triangulation where we study
the $\chi$-completeness hypothesis. Moreover, we show that $L_1$ and $L_2$ is not necessarily
essentially self-adjoint on the simple case.

We can extend the results in this paper to more general 2-simplicial complex, where the oriented faces
are not necessarily triangles. Particularly we can give a more general expression of the operator $d^1.$
More precisely, one can take 2-simplicial complexes with the number of edges of an oriented face bounded.
Indeed this hypothesis is important to give a meaning of the inequality in Definition \ref{defi4.2}.

\section{Preliminaries}
\subsection{The basic concepts}

A graph $\Kc$ is a pair ($\Vc$,$\Ec$), where $\Vc$
is the countable set of vertices and $\Ec$ the set of
oriented edges, considered as a subset of $\Vc\times\Vc$.
When two vertices $x$ and $y$ are connected by an edge $e$,
we say they are neighbors. We denote $x\sim y$ and $e=[x,y]\in\Ec.$
We assume that $\Ec$ is symmetric, \emph{ie.} $[x,y]\in\Ec\Rightarrow[y,x]\in\Ec.$
An oriented graph $\Kc$ is given by a partition of $\Ec:$
$$\Ec=\Ec^-\cup\Ec^+
$$
$$(x,y)\in\Ec^-\Leftrightarrow(y,x)\in\Ec^+
$$

In this case for $e=(x,y)\in\Ec^-,$ we define the origin $e^-=x,$
the termination $e^+=y$ and the opposite edge $-e=(y,x)\in\Ec^+.$
Let $c:\mathcal{V}\rightarrow(0,\infty)$ the weight on the vertices.
We also have $r:\Ec\longrightarrow (0,\infty)$ the weight on the oriented edges with
$$\forall e\in\Ec,~r(-e)=r(e).
$$

A \emph{path} between two vertices $x,y\in \Vc$ is a finite set of oriented edges
  $e_1, . . . , e_n, n\geq1$ such that
$$e^{-}_{1}=x,~e^{+}_{n}=y~\textrm{and},~\textrm{if}~n\geq2,
  ~\forall j,~1\leq j\leq n-1\Longrightarrow e^{+}_{j}=e^{-}_{j+1}.
$$

The path is called a \emph{cycle} or \emph{closed} when the origin
and the end are identical, \emph{ie.} $e^{-}_{1}=e^{+}_{n},$ with $n\geq3.$
If no cycles appear more than once in a path, the path is called a
\emph{simple path}. The graph $\Kc$ is \emph{connected} if any two vertices
$x$ and $y$ can be connected by a path with $e_1^-=x$ and $e_n^+=y.$
We say that the graph $\Kc$ is \emph{locally finite} if each vertex belongs
to a finite number of edges. The graph $\Kc$ is \emph{without loops} if there
is not the type of edges $(x,x),$ \emph{ie.}
$$\forall x\in\Vc \Longrightarrow (x,x)\notin \Ec.
$$
\subsubsection{} The set of \emph{neighbors} of $x\in \Vc$ is denoted by
$$\Vc(x):=\{y\in\Vc:y\sim x\}.
$$
\subsubsection{} The \emph{degree} of $x\in \Vc$ is by definition $deg(x)$,
the number of neighbors of $x$.
\subsubsection{} The \emph{combinatorial distance} $d_{comb}$ on $\Kc$ is
$$d_{comb}(x,y)=\min\{n,~\{e_{i}\}_{1\leq i\leq n}\subseteq\Ec
\mbox{ a path between the two vertices \emph{x} and \emph{y}}\}.
$$

\subsubsection{} Let $B$ be a finite subset of $\Vc.$ We define the \emph{edge boundary}
$\partial_{\Ec}B$ of $B$ by
$$\partial_{\Ec}B:=\{e\in\Ec\mbox{ such that }\{e^{-},e^{+}\}\cap B\neq \emptyset
\mbox{ and }\{e^{-},e^{+}\}\cap B^{c}\neq \emptyset\}.
$$

In the sequel, we assume that
\begin{center}
\textbf{$\Kc$ is without loops, connected, locally finite and oriented}
\end{center}
\begin{defi}\label{facebook}
An oriented face of $\Kc$ is a surface limited by a simple closed path,
considered as an element of $\Ec^n$ with $n\geqslant3,$ i.e
$$\varpi \mbox{ an oriented face }\Rightarrow\exists n\geq3,
~\varpi=(e_1,e_2,...,e_n)\in\Ec^n\mbox{ such that }\{e_{i}\}_{1\leq i\leq n}
\subseteq\Ec\mbox{ is a simple closed path}.
$$
\end{defi}

Let $\Fc$ be the set of all oriented faces of $\Kc,$ we consider
the pair ($\Kc$,$\Fc$) as a 2-simplicial complex, we denote it by $\Tc$.
We can denote also $\Tc=$($\Vc$,$\Ec$,$\Fc$).
\begin{rem}
Care should be taken not to confuse the simple cycles and the oriented faces.
Indeed, one can have simple cycles that are not oriented faces.
\end{rem}

For a face $\varpi=(e_1,e_2,...,e_n)\in\Fc,$ we have
$$\varpi=(e_i,...,e_n,e_1,...,e_{i-1})\in\Fc,~\forall 3\leq i\leq n-1.
$$

We can denote also
$$\varpi=(e_2,e_3,...,e_n,e_1)=....=(e_n,e_1,e_2,...,e_{n-1})\in\Fc.
$$

Because $\Kc$ is an oriented graph, we demand
$$(e_1,e_2,...,e_n)\in\Fc\Rightarrow(-e_n,-e_{n-1},...,-e_2,-e_1)\in\Fc.
$$

Given $\varpi=(e_1,e_2,...,e_n)\in\Fc,$ the \emph{opposite face}
of $\varpi$ is denoted by
$$-\varpi=(-e_n,-e_{n-1},...,-e_2,-e_1)\in\Fc.
$$


Let $B$ be a finite subset of $\Vc.$ We define the \emph{face boundary}
$\partial_{\Fc}B$ of $B$ by
$$\partial_{\Fc}B:=\{\sigma=(e_1,e_2,...,e_n)\in\Fc,~\exists i
\mbox{ such that }e_i\in \partial_{\Ec}B,~n\geqslant 3\}.
$$

\begin{defi}\emph{(Triangulation)}
\emph{A triangulation} is a 2-simplicial complex such that all the faces are triangles.
\end{defi}

\begin{rem}
In the definition of a triagulation we demand that faces are triangles.
In the sequel, we assume also that each triangular cycle is an oriented face for simplicity reasons.
Indeed all the results of this work can be extended easily to any triangulation.
\end{rem}

In the sequel we will represent the oriented faces by their vertices
$$\varpi=(e_1,e_2,e_3)=[e^-_1=e^+_3,e^+_1=e^-_2,e^+_2=e^-_3]\in\Fc.
$$

For a face $\varpi=[x,y,z]\in\Fc.$ Let us set
$$\varpi=[x,y,z]=[y,z,x]=[z,x,y]\in\Fc\Rightarrow -\varpi=[y,x,z]=[x,z,y]=[z,y,x]\in\Fc.
$$

To define weighted triangulations we need weights, let us give $s:\Fc\rightarrow(0,\infty)$ the weight on oriented faces
such that for all $\varpi\in\Fc,~s(-\varpi)=s(\varpi).$ The weighted triangulation $(\Tc,c,r,s)$ is given by the triangulation
$\Tc=(\Vc,\Ec,\Fc).$ We say that $\Tc$ is \emph{simple} if the weights of the vertices, the edges and faces equals $1.$
For an edge $e\in\Ec$, we also denote the oriented face $[e^-,e^+,x]$ by $(e,x),$ with $x\in\Vc(e^-)\cap\Vc(e^+).$
The set of vertices belonging to the edge $e\in\Ec$ is given by
$$\Fc_e:=\{x\in\Vc,~(e,x)\in\Fc\}=\Vc(e^-)\cap\Vc(e^+).
$$

\subsection{Functions spaces} We denote the set of 0-cochains or functions on $\Vc$ by:
$$\Cc(\Vc)=\{f:\Vc\rightarrow\C\}
$$
and the set of functions of finite support by $\Cc_{c}(\Vc).$

Similarly, we denote the set of 1-cochains or 1-forms on $\Ec$ by:
$$\Cc(\Ec)=\{\varphi:\Ec\rightarrow\C,~\varphi(-e)=-\varphi(e)\}
$$
and the set of 1-forms of finite support by $\Cc_{c}(\Ec).$

Moreover, we denote the set of 2-cochains or 2-forms on $\Fc$ by:
$$\Cc(\Fc)=\{\phi:\Fc\rightarrow\C,~\phi(-\varpi)=-\phi(\varpi)\}
$$
and the set of 2-forms of finite support by $\Cc_{c}(\Fc).$

Let us define the Hilbert spaces $l^2(\Vc),~l^{2}(\Ec)$ and $l^{2}(\Fc)$ as the sets of
cochains with finite norm, we have
\begin{itemize}
\item[(a)]
$$l^2(\Vc):=\{f\in\Cc(\Vc);~\displaystyle\sum_{x\in\Vc}c(x)|f(x)|^{2}<\infty\},
$$
with the inner product
$$\langle f,g \rangle_{l^2(\Vc)}:=\displaystyle\sum_{x\in\Vc}c(x)f(x)\overline{g}(x).
$$
\item[(b)]
$$l^{2}(\Ec):=\{\varphi\in\Cc(\Ec);~\displaystyle\sum_{e\in\Ec}r(e)|\varphi(e)|^{2}<\infty\},
$$
with the inner product
$$\langle \varphi,\psi \rangle_{l^{2}(\Ec)}:=\frac{1}{2}
\displaystyle\sum_{e\in\Ec}r(e)\varphi(e)\overline{\psi}(e).
$$
\item[(c)]
$$l^{2}(\Fc):=\{\phi\in\Cc(\Fc);~\displaystyle\sum_{\varpi\in\Fc}s(\varpi)|\phi(\varpi)|^{2}<\infty\},
$$
with the inner product
$$\langle \phi_{1},\phi_{2} \rangle_{l^{2}(\Fc)}
=\dfrac{1}{6}\displaystyle\sum_{[x,y,z]\in\Fc}s(x,y,z)\phi_{1}(x,y,z)
\overline{\phi_{2}}(x,y,z).
$$
\end{itemize}



The direct sum of the spaces $l^{2}(\Vc)$, $l^{2}(\Ec)$ and $l^{2}(\Fc)$
can be considered as a new Hilbert space denoted by $\Hc$, that is
$$\Hc=l^{2}(\Vc)\oplus l^{2}(\Ec)\oplus l^{2}(\Fc),
$$
with the norm
$$\forall F=(f,\varphi,\phi)\in\mathcal{H},~\|F\|_{\Hc}^{2}=\|f\|_{l^{2}(\Vc)}^{2}
+\|\varphi\|_{l^{2}(\Ec)}^{2}+\|\phi\|_{l^{2}(\Fc)}^{2}.
$$

\subsection{Operators}
We give in this part the expressions of the operators introduced
on graphs which are already well known and we also give other operators
acting on triangulations.
\subsubsection{The difference operator}

By analogy to electric networks of voltage differences across edges
leading to currents \cite{LP}, we define \emph{the difference operator}
$d^{0}:\Cc_{c}(\Vc)\longrightarrow \Cc_{c}(\Ec)$ by
  $$\forall f\in \Cc_{c}(\Vc),~d^{0}(f)(e)=f(e^{+})-f(e^{-}).
  $$
\subsubsection{The co-boundary operator}
It is the formal adjoint of $d^0,$  denoted
$\delta^0:\Cc_{c}(\Ec)\longrightarrow \Cc_{c}(\Vc),$ (see \cite{AT})
acts as
$$\forall \varphi\in \Cc_{c}(\Ec),~\delta^{0}(\varphi)(x)=\dfrac{1}{c(x)}
\displaystyle\sum_{e,e^{+}=x}r(e)\varphi(e).
$$

\subsubsection{The exterior derivative}
It is the operator $d^{1}:\Cc_{c}(\Ec)\longrightarrow \Cc_{c}(\Fc),$ given by
$$\forall \psi\in \Cc_{c}(\Ec),~d^{1}(\psi)(x,y,z)=\psi(x,y)+\psi(y,z)+\psi(z,x).
$$

\subsubsection{The co-exterior derivative} It is the formal adjoint of $d^1,$  denoted
$\delta^1:\Cc_{c}(\Fc)\longrightarrow \Cc_{c}(\Ec),$
which satisfies
\begin{align} \label{Eq2.1}
\langle d^{1}\psi,\phi\rangle_{l^{2}(\Fc)}
=\langle \psi,\delta^{1}\phi\rangle_{l^{2}(\Ec)}
,~\forall (\psi,\phi)\in\Cc_{c}(\Ec)\times\Cc_{c}(\Fc).
\end{align}
\begin{lem} \label{lem1}
The formal adjoint $\delta^{1}:\Cc_{c}(\Fc)\longrightarrow \Cc_{c}(\Ec)$, is given by
$$\delta^{1}(\phi)(e)
=\dfrac{1}{r(e)}\sum_{x\in \Fc_{e}}s(e,x)\phi(e,x).
$$
\end{lem}
\begin{dem}

Let $(\psi,\phi)\in \Cc_{c}(\Ec)\times\Cc_{c}(\Fc).$ The equation (\ref{Eq2.1}) gives
\begin{equation*}
\begin{split}
\langle d^1\psi,\phi \rangle_{l^2(\Fc)}
& =\dfrac{1}{6}\displaystyle\sum_{[x,y,z]\in\Fc}s(x,y,z)d^1(\psi)(x,y,z)\overline{\phi}(x,y,z)\\
& =\dfrac{1}{2}\displaystyle\sum_{[x,y,z]\in\Fc}s(x,y,z)\psi(x,y)\overline{\phi}(x,y,z)\\
& =\langle \psi,\delta^1\phi \rangle_{l^2(\Ec)}.
\end{split}
\end{equation*}

To justify it note that the expression of $d^1$ contributing to the first sum is divided into
three similar parts. So it remains to show only
\begin{equation*}
\begin{split}
\displaystyle\sum_{[x,y,z]\in\Fc}s(x,y,z)\psi(x,y)\overline{\phi}(x,y,z)
& =\displaystyle\sum_{e\in\Ec}\psi(e)\displaystyle\sum_{x\in\Fc_{e}}
s(e,x)\overline{\phi}(e,x)\\
& =\displaystyle\sum_{e\in\Ec}r(e)\psi(e)
\overline{\left(\dfrac{1}{r(e)}\sum_{x\in\Fc_{e}}s(e,x)\phi(e,x)\right)}
\end{split}
\end{equation*}
\end{dem}
\subsubsection{Gau\ss-Bonnet operator on $\Tc$}
By analogy to Riemannian geometry, we use the decomposition of the operators in \cite{EL}
to define the Gau\mbox{\ss}-Bonnet operator. Let us begin by defining the operator
$$d:\Cc_{c}(\Vc)\oplus\Cc_{c}(\Ec)\oplus\Cc_{c}(\Fc)\circlearrowleft
$$
by
$$\forall(f,\varphi,\phi)\in\Cc_{c}(\Vc)\oplus\Cc_{c}(\Ec)\oplus\Cc_{c}(\Fc),
~d(f,\varphi,\phi)=(0,d^{0}f,d^{1}\varphi),
$$
and $\delta$ the formal adjoint of $d.$ Thus it satisfies
\begin{align} \label{Eq2.2}
\langle d(f_{1},\varphi_{1},\phi_{1}),(f_{2},\varphi_{2},\phi_{2})\rangle_{\Hc}
=\langle(f_{1},\varphi_{1},\phi_{1}),\delta(f_{2},\varphi_{2},\phi_{2})\rangle_{\Hc},
\end{align}
for all $(f_{1},\varphi_{1},\phi_{1}),(f_{2},\varphi_{2},\phi_{2})\in\Cc_{c}(\Vc)
\oplus\Cc_{c}(\Ec)\oplus\Cc_{c}(\Fc)$.
\begin{lem}
Let $\Tc=(\Kc,\Fc)$ be a triangulation. Then
$$\delta:\Cc_{c}(\Vc)\oplus\Cc_{c}(\Ec)\oplus\Cc_{c}(\Fc)\circlearrowleft
$$
is given by
$$\delta(f,\varphi,\phi)=(\delta^0\varphi,\delta^{1}\phi,0),
~\forall(f,\varphi,\phi)\in\Cc_{c}(\Vc)\oplus\Cc_{c}(\Ec)\oplus\Cc_{c}(\Fc).
$$
\end{lem}
\begin{dem}
Let $(f_{1},\varphi_{1},\phi_{1}),(f_{2},\varphi_{2},\phi_{2})\in \Cc_{c}(\Vc)
\oplus\Cc_{c}(\Ec)\oplus\Cc_{c}(\Fc)$. Using the equation (\ref{Eq2.2})
\begin{equation*}
\begin{split}
\langle d(f_{1},\varphi_{1},\phi_{1}),(f_{2},\varphi_{2},\phi_{2})\rangle_{\Hc}
& =\langle (0,d^{0}f_{1},d^{1}\varphi_{1}),(f_{2},\varphi_{2},\phi_{2})\rangle_{\Hc}\\
& =\langle d^{0}f_{1},\varphi_{2}\rangle_{l^{2}(\Ec)}
+\langle d^{1}\varphi_{1},\phi_{2}\rangle_{l^{2}(\Fc)}\\
& =\langle f_{1},\delta^0\varphi_{2}\rangle_{l^{2}(\Vc)}
+\langle \varphi_{1},\delta^{1}\phi_{2}\rangle_{l^{2}(\Ec)}\\
&=\langle(f_{1},\varphi_{1},\phi_{1}),(\delta^0\varphi_{2},\delta^{1}\phi_{2},0)\rangle_{\Hc}.
\end{split}
\end{equation*}
\end{dem}
\begin{defi}
Let $\Tc=(\Kc,\Fc)$ be a triangulation, the Gau\ss -Bonnet operator defined as
$$T:=d+\delta:\Cc_{c}(\Vc)\oplus\Cc_{c}(\Ec)\oplus\Cc_{c}(\Fc)\circlearrowleft
$$
is given by
$$T(f,\varphi,\phi)=(\delta^0\varphi,d^{0}f+\delta^{1}\phi,d^{1}\varphi)
$$
for all $(f,\varphi,\phi)\in\Cc_{c}(\Vc)\oplus\Cc_{c}(\Ec)\oplus\Cc_{c}(\Fc).$
Moreover, the matrix representation of $T$ is given by
\begin{center}
$T\equiv\begin{pmatrix}
0&\delta^0&0 \\
d^{0}&0&\delta^{1} \\
0&d^{1}&0
\end{pmatrix}$
\end{center}
\end{defi}

\begin{lem} \label{lem2}
If $\Tc=(\Kc,\Fc)$ is a triangulation then $d^{1}d^{0}=\delta^0\delta^{1}=0.$
\end{lem}
\begin{dem}
Let $f\in\Cc_{c}(\Vc),$ we have that
\begin{equation*}
\begin{split}
d^{1}(d^{0}f)(x,y,z)
& =d^{0}f(x,y)+d^{0}f(y,z)+d^{0}f(z,x)\\
& =(f(y)-f(x))+(f(z)-f(y))+(f(x)-f(z))=0.
\end{split}
\end{equation*}
Since $d^{1}d^{0}=0$ and the operator $\delta^0\delta^{1}$ is the formal adjoint of
$d^{1}d^{0}.$ Then $\delta^0\delta^{1}=0.$
\end{dem}

Before giving an important result for $f\in \Cc(\Vc),$ we define the two operators
$~\widetilde{}:\Cc(\Vc)\longrightarrow\Cc(\Ec)$ by $f\mapsto\widetilde{f}$ and
$~\widetilde{\widetilde{}}:\Cc(\Vc)\longrightarrow\Cc(\Fc)$ by $f\mapsto\widetilde{\widetilde{f}},$ where
$$\widetilde{f}(e):=\frac{1}{2}\left(f(e^+)+f(e^-)\right).
$$
$$\widetilde{\widetilde{f}}(x,y,z):=\frac{1}{3}\left(\widetilde{f}(x,y)
+\widetilde{f}(y,z)+\widetilde{f}(z,x)\right)
=\frac{1}{3}(f(x)+f(y)+f(z)).
$$

\emph{The exterior product} of two 1-forms defined as $.\wedge_{disc}.:\Cc(\Ec)\times\Cc(\Ec)\longrightarrow\Cc(\Fc),$
is given by:
\begin{equation*}
\begin{split}
\left(\psi\wedge_{disc}\varphi\right)(x,y,z)&=
\left[\psi(z,x)+\psi(z,y)\right]\varphi(x,y)\\
&+\left[\psi(x,y)+\psi(x,z)\right]\varphi(y,z)\\
&+\left[\psi(y,z)+\psi(y,x)\right]\varphi(z,x).
\end{split}
\end{equation*}

It satisfies $\psi\wedge_{disc}\varphi=-\left(\varphi\wedge_{disc}\psi\right)
=-\varphi\wedge_{disc}\psi=\varphi\wedge_{disc}-\psi,$ for all $\varphi,\psi\in\Cc(\Ec).$
\begin{lem}\emph{(Derivation properties)}\label{lem3}
Let $(f,\varphi,\phi)\in\Cc_{c}(\Vc)\times\Cc_{c}(\Ec)\times\Cc_{c}(\Fc).$ Then
\begin{equation}\label{Eq4.1}
d^1(\widetilde{f}\varphi)(x,y,z)=\widetilde{\widetilde{f}}(x,y,z)d^1(\varphi)(x,y,z)
+\frac{1}{6}\left(d^0(f)\wedge_{disc}\varphi\right)(x,y,z).
\end{equation}
\begin{equation}\label{Eq4.2}
\delta^1(\widetilde{\widetilde{f}}\phi)(e)=\widetilde{f}(e)\delta^1(\phi)(e)
+\dfrac{1}{6r(e)}\displaystyle\sum_{x\in\Fc_e}s(e,x)
\left[d^0(f)(e^-,x)+d^0(f)(e^+,x)\right]\phi(e,x).
\end{equation}
\end{lem}
\begin{dem}
\begin{enumerate}
  \item [(1)]Let $(f,\varphi)\in\Cc_{c}(\Vc)\times\Cc_{c}(\Ec),$ we have
\begin{equation*}
\begin{split}
d^1(\widetilde{f}\varphi)(x,y,z)
& = \widetilde{f}(x,y)\varphi(x,y)+\widetilde{f}(y,z)\varphi(y,z)
+\widetilde{f}(z,x)\varphi(z,x)\\
& = [\widetilde{f}(x,y)+\widetilde{f}(y,z)+\widetilde{f}(z,x)]
[\varphi(x,y)+\varphi(y,z)+\varphi(z,x)]\\
& - \left(\widetilde{f}(y,z)+\widetilde{f}(z,x)\right)\varphi(x,y)-
\left(\widetilde{f}(z,x)+\widetilde{f}(x,y)\right)\varphi(y,z)\\
& - \left(\widetilde{f}(x,y)+\widetilde{f}(y,z)\right)\varphi(z,x)\\
& = \widetilde{\widetilde{f}}(x,y,z)d^1(\varphi)(x,y,z)\\
& + \left(\frac{2}{3}\widetilde{f}(x,y)-\frac{1}{3}
\left[\widetilde{f}(y,z)+\widetilde{f}(z,x)\right]\right)\varphi(x,y)\\
& + \left(\frac{2}{3}\widetilde{f}(y,z)-\frac{1}{3}
\left[\widetilde{f}(z,x)+\widetilde{f}(x,y)\right]\right)\varphi(y,z)\\
& + \left(\frac{2}{3}\widetilde{f}(z,x)-\frac{1}{3}
\left[\widetilde{f}(x,y)+\widetilde{f}(y,z)\right]\right)\varphi(z,x).
\end{split}
\end{equation*}

On the other hand, we have
\begin{equation*}
\begin{split}
\left(\frac{2}{3}\widetilde{f}(x,y)-\frac{1}{3}
\left[\widetilde{f}(y,z)+\widetilde{f}(z,x)\right]\right)
& = \frac{1}{3}\left(\left[\widetilde{f}(x,y)-\widetilde{f}(y,z)\right]
+\left[\widetilde{f}(x,y)-\widetilde{f}(z,x)\right]\right)\\
& = \frac{1}{6}\left(d^{0}(f)(z,x)
+d^{0}(f)(z,y)\right).
\end{split}
\end{equation*}

Similarly, we get
$$\left(\frac{2}{3}\widetilde{f}(y,z)-\frac{1}{3}
\left[\widetilde{f}(z,x)+\widetilde{f}(x,y)\right]\right)\varphi(y,z)
=\frac{1}{6}\left(d^{0}(f)(x,y)
+d^{0}(f)(x,z)\right)\varphi(y,z).
$$
and
$$\left(\frac{2}{3}\widetilde{f}(z,x)-\frac{1}{3}
\left[\widetilde{f}(x,y)+\widetilde{f}(y,z)\right]\right)\varphi(z,x)
=\frac{1}{6}\left(d^{0}(f)(y,z)
+d^{0}(f)(y,x)\right)\varphi(z,x).
$$

Hence, we have
\begin{equation*}
\begin{split}
d^1(\widetilde{f}\varphi)(x,y,z)
& = \widetilde{\widetilde{f}}(x,y,z)d^1(\varphi)(x,y,z)\\
&+\frac{1}{6}\left[d^{0}(f)(z,x)+d^{0}(f)(z,y)\right]\varphi(x,y)\\
&+\frac{1}{6}\left[d^{0}(f)(x,y)+d^{0}(f)(x,z)\right]\varphi(y,z)\\
&+\frac{1}{6}\left[d^{0}(f)(y,z)+d^{0}(f)(y,x)\right]\varphi(z,x).
\end{split}
\end{equation*}
  \item [(2)]Let $(f,\phi)\in\Cc_{c}(\Vc)\times\Cc_{c}(\Fc).$ Then by Lemma \ref{lem1},
  \begin{equation*}
  \begin{split}
  \delta^1(\widetilde{\widetilde{f}}\phi)(e)
  & = \dfrac{1}{r(e)}\displaystyle\sum_{x\in\Fc_e}s(e,x)
  \widetilde{\widetilde{f}}(e,x)\phi(e,x)\\
  & = \dfrac{1}{3r(e)}\displaystyle\sum_{x\in\Fc_e}s(e,x)
  \left(\widetilde{f}(e)+\widetilde{f}(e^+,x)+\widetilde{f}(x,e^-)\right)
  \phi(e,x)\\
  & = \frac{1}{3}\widetilde{f}(e)\delta^1(\phi)(e)
    + \dfrac{1}{3r(e)}
      \displaystyle\sum_{x\in\Fc_e}s(e,x)
      \widetilde{f}(e^+,x)\phi(e,x)\\
  & + \dfrac{1}{3r(e)}\displaystyle\sum_{x\in\Fc_e}s(e,x)
      \widetilde{f}(x,e^-)\phi(e,x).\\
  & = \widetilde{f}(e)\delta^1(\phi)(e)\\
  & + \dfrac{1}{3r(e)}
      \displaystyle\sum_{x\in\Fc_e}s(e,x)
      \left(\widetilde{f}(e^+,x)-\widetilde{f}(e)\right)\phi(e,x)\\
  & + \dfrac{1}{3r(e)}
      \displaystyle\sum_{x\in\Fc_e}s(e,x)
      \left(\widetilde{f}(x,e^-)-\widetilde{f}(e)\right)\phi(e,x)\\
  & = \widetilde{f}(e)\delta^1(\phi)(e)\\
  & + \dfrac{1}{6r(e)}
      \displaystyle\sum_{x\in\Fc_e}s(e,x)
      \left[d^0(f)(e^-,x)+d^0(f)(e^+,x)\right]\phi(e,x).
  \end{split}
  \end{equation*}
\end{enumerate}
\end{dem}

\subsubsection{Laplacian}
Through the Gau\ss -Bonnet operator $T,$ we can define the discrete Laplacian on $\Tc$.
So, Lemma \ref{lem2} induces the following definition
\begin{defi}
Let $\Tc=(\Kc,\Fc)$ be a triangulation, the Laplacian on $\Tc$ defined as
$$L:=T^{2}:\Cc_{c}(\Vc)\oplus\Cc_{c}(\Ec)\oplus\Cc_{c}(\Fc)\circlearrowleft
$$
is given by
$$L(f,\varphi,\phi)=(\delta^0d^{0}f,(d^{0}\delta^0+\delta^{1}d^{1})\varphi,d^{1}\delta^{1}\phi).
$$
for all $(f,\varphi,\phi)\in\Cc_{c}(\Vc)\oplus\Cc_{c}(\Ec)\oplus\Cc_{c}(\Fc).$
\end{defi}
\begin{rem}
We can write
  $$L:=L_{0}\oplus L_{1}\oplus L_{2},
  $$
where $L_{0}$ is the discrete Laplacian acting on functions given by
$$L_{0}(f)(x):=\delta^0d^{0}(f)(x)=\dfrac{1}{c(x)}
\displaystyle\sum_{e,e^{+}=x}r(e)d^{0}(f)(e),
$$
with $f\in\Cc_{c}(\Vc),$ and where $L_{1}$ is the discrete Laplacian
acting on 1-forms given by
\begin{equation*}
\begin{split}
L_{1}(\varphi)(x,y)
& :=(d^{0}\delta^0+\delta^{1}d^{1})(\varphi)(x,y)\\
& =\dfrac{1}{c(y)}\displaystyle\sum_{e,e^{+}=y}r(e)\varphi(e)
-\dfrac{1}{c(x)}\displaystyle\sum_{e,e^{+}=x}r(e)\varphi(e)
+\dfrac{1}{r(x,y)}\displaystyle\sum_{z\in\Fc_{[x,y]}}s(x,y,z)d^{1}(\varphi)(x,y,z),
\end{split}
\end{equation*}
with $\varphi\in\Cc_{c}(\Ec),$ and where also $L_{2}$ is the discrete
Laplacian acting on 2-forms given by
\begin{equation*}
\begin{split}
L_{2}(\phi)(x,y,z)
&:=d^{1}\delta^{1}(\phi)(x,y,z)\\
& =\dfrac{1}{r(x,y)}\displaystyle\sum_{u\in\Fc_{[x,y]}}s(x,y,u)\phi(x,y,u)\\
& +\dfrac{1}{r(y,z)}\displaystyle\sum_{u\in\Fc_{[y,z]}}s(y,z,u)\phi(y,z,u)\\
& +\dfrac{1}{r(z,x)}\displaystyle\sum_{u\in\Fc_{[z,x]}}s(z,x,u)\phi(z,x,u),
\end{split}
\end{equation*}
with $\phi\in\Cc_{c}(\Fc).$
\end{rem}
\begin{rem}
The operator $L_1$ is called the \emph{full Laplacian} and defined as $L_1=L_1^-+L_1^+,$
where $L_1^-=d^{0}\delta^0$ \emph{(}resp. $L_{1}^{+}=\delta^{1}d^{1}$\emph{)} is called
the \emph{lower Laplacian} \emph{(}resp. the \emph{upper Laplacian}\emph{)}. In both articles
\cite{AT} and \cite{BGJ}, the authors denote $\Delta_0=\delta^0d^0$ and $\Delta_1=d^0\delta^0.$
In this work, we have $L_0=\Delta_0$ and $L_{1}^{-}=\Delta_1.$
\end{rem}
\section{Closability}\label{Section3}
On a connected locally finite graph, the operators $d^{0}$ and $\delta^0$
are closable (see \cite{AT}). The next Lemma proves the closability of the
operators $d^{1}$ and $\delta^{1}$ on a triangulation.
\begin{lem}
Let $\Tc=(\Kc,\Fc)$ be a weighted triangulation. Then the operators $d^{1}$ and $\delta^{1}$ are closable.
\end{lem}
\begin{dem}
\begin{itemize}
\item Let $(\varphi_{n})_{n\in\N}$ be a sequence from $\Cc_{c}(\Ec)$ and $\phi\in l^{2}(\Fc)$ such that
$$\displaystyle\lim_{n\rightarrow\infty}\left(\|\varphi_{n}\|_{l^{2}(\Ec)}+\|d^{1}\varphi_{n}-\phi\|_{l^{2}(\Fc)}\right)=0,
$$
then for each edge $e,~\varphi_{n}(e)$ converges to $0$ and for each face $\varpi,~d^{1}(\varphi_{n})(\varpi)$
converges to $\phi(\varpi).$ But by the expression of $d^{1}$ and local finiteness of $\Tc,$ for each face
$\varpi,~d^{1}(\varphi_{n})(\varpi)$ converges to $0.$ Thus we have that $\phi=0.$
\item The same can be done for $\delta^{1}:$ Let $(\phi_{n})_{n\in\N}$ be a sequence from $\Cc_{c}(\Fc)$ and
$\varphi\in l^{2}(\Ec)$ such that
$$\displaystyle\lim_{n\rightarrow\infty}\left(\|\phi_{n}\|_{l^{2}(\Fc)}+\|\delta^{1}\phi_{n}-\varphi\|_{l^{2}(\Ec)}\right)=0,
$$
then for each face $\sigma,~\phi_{n}(\sigma)$ converges to $0$ and for each edge $e$, $\delta^{1}(\phi_{n})(e)$
converges to $\varphi(e).$ But by the expression of $\delta^{1}$ and local finiteness of $\Tc$, for each edge $e$,
$\delta^{1}(\phi_{n})(e)$ converges to $0.$ Thus we have that $\varphi=0.$
\end{itemize}
\end{dem}

The smallest extension is the closure (see \cite{S},\cite{RSv1}), denoted
$\overline{d^{0}}:=d^{0}_{min}$ (resp. $\overline{\delta^0}:=\delta_{min}^0$,
$\overline{d^{1}}:=d^{1}_{min}$, $\overline{\delta^{1}}:=\delta^{1}_{min},$
$\overline{T}:=T_{min}$, $\overline{L}:=L_{min}$) has the domain
$$Dom(d^{0}_{min})=\left\{f\in l^{2}(\Vc);~\exists (f_{n})_{n\in\N},
~f_{n}\in\Cc_{c}(\Vc),
~\displaystyle\lim_{n\rightarrow\infty}\|f_{n}-f\|_{l^{2}(\Vc)}=0,
~\displaystyle\lim_{n\rightarrow\infty}d^{0}(f_{n})\mbox{ exists in }l^{2}(\Ec)\right\},
$$
for such an $f$, one puts
$$d^{0}_{min}(f)=\displaystyle\lim_{n\rightarrow\infty}d^{0}(f_{n}).
$$

We notice that $d^0_{min}(f)$ is independent of the sequence $(f_n)_{n\in\N},$ because $d^{0}$ is closable.

The largest is $d^{0}_{max}=(\delta^0)^{*}$, the adjoint operator of
$\delta^0_{min}$, (resp. $\delta^0_{max}=(d^0)^{*}$, the adjoint operator of $d^{0}_{min}$).

We also note $d^{1}_{max}=(\delta^1)^{*}$, the adjoint operator of $\delta^{1}_{min}$,
(resp $\delta^{1}_{max}=(d^1)^{*}$, the adjoint operator of $d^{1}_{min}$).
\begin{propr} \label{pro1}
Let $\Tc=(\Kc$,$\Fc)$ be a weighted triangulation. Then
$$Dom(T_{min})\subseteq Dom(d^{0}_{min})\oplus\left(Dom(\delta^{0}_{min})\cap Dom(d^{1}_{min})\right)
\oplus Dom(\delta^{1}_{min}).
$$
\end{propr}
\begin{dem}

Let $F=(f,\varphi,\phi)\in Dom(T_{min}),$ so
there exists a sequence
$(F_n)_n=\left((f_n,\varphi_n,\phi_n)\right)_n\subseteq \Cc_{c}(\Vc)
\oplus\Cc_{c}(\Ec)\oplus\Cc_{c}(\Fc)$ such that
$\displaystyle\lim_{n\rightarrow\infty}F_n=F$ in $\Hc$
and $(TF_n)_{n\in\N}$ converges in $\Hc.$ Let us denote by
$l_0=(f_0,\varphi_0,\phi_0)$ this limit. Therefore
\begin{equation*}
\begin{split}
\|TF_n-l_0\|_{\Hc}^2
& =\|\delta^0\varphi_n-f_0\|_{l^2(\Vc)}^2
  +\|(d^0+\delta^1)(f_n,\phi_n)-\varphi_0\|_{l^2(\Ec)}^2
  +\|d^1\varphi_n-\phi_0\|_{l^2(\Fc)}^2.
\end{split}
\end{equation*}

Hence $\delta^0\varphi_n\rightarrow f_0$ and $d^1\varphi_n\rightarrow\phi_0$
respectively in $l^2(\Vc)$ and in $l^2(\Ec).$
So, by definition, $\varphi\in Dom(\delta^{0}_{min})\cap Dom(d^{1}_{min}),$
$f_0=\delta_{min}^0\varphi$ and $\phi_0=d_{min}^1\varphi.$
Moreover, we combine the parallelogram identity with Lemma \ref{lem2}
to obtain the following result
$$\|(d^0+\delta^1)(f,\phi)\|_{l^2(\Ec)}^2=\|d^0(f)\|_{l^2(\Ec)}^2
+\|\delta^1(\phi)\|_{l^2(\Ec)}^2,~\forall n\in\N,~ \forall(f,\phi)\in\Cc_{c}(\Vc)\times\Cc_{c}(\Fc).
$$

Since $\left((d^0+\delta^1)(f_n,\phi_n)\right)_n$ converges in $l^2(\Ec),$
then by completeness of $l^2(\Ec)$ $\left(d^0(f_n)\right)_n$ and $\left(\delta^1(\phi_n)\right)_n$
are convergent in $l^2(\Ec).$ Thus, we conclude that $f\in Dom(d^{0}_{min})$ and $\phi\in Dom(\delta^{1}_{min}).$
\end{dem}
\begin{propr}\label{Yasch}
Let $\Tc=(\Kc,\Fc)$ be a weighted triangulation. Then
$$Dom(L_{min})\subseteq Dom(\delta_{min}^0d_{min}^0)\oplus \left(Dom(d_{min}^0\delta_{min}^0)\cap Dom(\delta_{min}^1d_{min}^1)\right)
\oplus Dom(d_{min}^1\delta_{min}^1).
$$
\end{propr}
\begin{dem}
\begin{enumerate}
\item[i)]We will show that $(L_0)_{min}\subseteq\delta^{0}_{min}d^{0}_{min}.$
First, we note that
$$Dom(\delta^{0}_{min}d^{0}_{min})=\{f\in Dom(d^{0}_{min}),
~d^{0}_{min}f\in Dom(\delta^{0}_{min})\}.
$$

Let $f\in Dom((L_0)_{min}),$ so there exists a sequence
$(f_n)_n\subseteq\Cc_{c}(\Vc)$ such that
$$f_n\rightarrow f\mbox{ in }l^2(\Vc),
~\delta^{0}d^{0}f_n\rightarrow (\delta^{0}d^{0})_{min}f \mbox{ in }l^2(\Vc).
$$

So, $(L_0f_n)_n$ is a Cauchy sequence. Moreover, we have
\begin{equation*}
\begin{split}
\forall n,m\in\N,~\|d^{0}f_n-d^{0}f_m\|_{l^2(\Ec)}^2
&=\langle d^{0}(f_n-f_m),d^{0}(f_n-f_m)\rangle_{l^{2}(\Ec)}\\
&=\langle\delta^{0}d^{0}(f_n-f_m),f_n-f_m \rangle_{l^{2}(\Vc)}\\
&=\langle L_{0}(f_n-f_m),f_n-f_m \rangle_{l^{2}(\Vc)}\\
&\leq\|L_{0}(f_n-f_m)\|_{l^2(\Vc)}\|f_n-f_m\|_{l^2(\Vc)}.
\end{split}
\end{equation*}

Thus $(d^{0}f_n)_n$ is a Cauchy sequence because $(L_{0}f_{n})_{n}$
is a Cauchy sequence and $(f_{n})_{n}$ is convergent. So, it is convergent in $l^2(\Ec).$
By closability, we conclude that $f\in Dom(\delta^{0}_{min}d^{0}_{min}).$
\item[ii)]First, for all $\varphi,\psi\in \Cc_{c}(\Ec)$ we have
\begin{equation}\label{*}
\langle L_1\varphi,\psi \rangle_{l^2(\Ec)}
=\langle (L_1^{-}+L_1^{+})\varphi,\psi \rangle_{l^2(\Ec)}
=\langle \delta^0\varphi,\delta^0\psi \rangle_{l^2(\Vc)}
+\langle d^1\varphi,d^1\psi \rangle_{l^2(\Fc)}.
\end{equation}

Using the same method as in i) with (\ref{*}) we obtain that
$(L_1^{-})_{min}\subseteq d^{0}_{min}\delta^{0}_{min}$ and
$(L_1^{+})_{min}\subseteq \delta^{1}_{min}d^{1}_{min}.$
It remains to show that we have
$$(L_1)_{min}\subseteq (L_1^{-})_{min}+(L_1^{+})_{min}
$$

Let $\varphi\in Dom((L_1)_{min}),$ so there exists a sequence
$(\varphi_n)_n\subseteq\Cc_{c}(\Ec)$ such that
$\varphi=\displaystyle\lim_{n\rightarrow\infty}\varphi_n\mbox{ in }l^2(\Ec)$
and $(L_1\varphi_n)_{n\in\N}$ converges in $l^2(\Ec).$
Then, by the parallelogram identity with Lemma \ref{lem2} we obtain
$$\|(L_1^{-}+L_1^{+})(\varphi_n)\|_{l^2(\Ec)}^2
=\|L_1^{-}(\varphi_n)\|_{l^2(\Ec)}^2
+\|L_1^{+}(\varphi_n)\|_{l^2(\Ec)}^2,~\forall n\in\N.
$$

Then $\left(L_1^{-}(\varphi_n)\right)_n$ and $\left(L_1^{+}(\varphi_n)\right)_n$
are convergent in $l^2(\Ec).$ Moreover, by the closability of
$L_1^{-}$ and $L_1^{+},$ we conclude that
$\varphi\in Dom((L_1^{-})_{min})\cap Dom((L_1^{+})_{min}).$
\item[iii)]Since, for any $\phi,\Theta\in\Cc_{c}(\Fc),$ we have
\begin{equation}\label{*e}
\langle L_2\phi,\Theta \rangle_{l^2(\Fc)}
=\langle d^{1}\delta^{1}\phi,\Theta \rangle_{l^2(\Fc)}
=\langle \delta^1\phi,\delta^1\Theta \rangle_{l^2(\Ec)}.
\end{equation}

Using the same method as in i) with (\ref{*e}) we obtain that
$(L_2)_{min}\subseteq d^{1}_{min}\delta^{1}_{min}.$
\end{enumerate}
\end{dem}

\section{Geometric hypothesis}
\subsection{$\chi$-completeness}

In this subsection, we give the geometric hypothesis for the triangulation $\Tc.$
First we recall the definition of $\chi$-completeness given in \cite{AT} for the case of graphs.
A graph $\Kc=(\Vc,\Ec)$ is $\chi$-complete if there exists an increasing sequence
of finite sets $(B_{n})_{n\in\N}$ such that $\Vc=\displaystyle\cup_{n\in\N}B_{n}$
and there exist related functions $\chi_{n}$ satisfying the following three conditions:
\begin{enumerate}
\item[i)] $\chi_{n}\in \Cc_{c}(\Vc),~0\leq\chi_{n}\leq1.$
\item[ii)] $x\in B_{n}\Rightarrow \chi_{n}(x)=1.$
\item[iii)] $\exists C>0$ such that $\forall n\in \N,~x\in \Vc$
$$\dfrac{1}{c(x)}\displaystyle\sum_{e\in\Ec,e^{\pm}=x}r(e)
|d^{0}\chi_{n}(e)|^{2}\leq C.
$$
\end{enumerate}
\begin{rem}
The $\chi$-completeness is related to the notion of intrinsic metric for weighted graphs.
This geometric hypothesis covers many situations that have been already studied.
Particularly in \cite{AT}, the authors prove that it is satisfied by locally finite graphs which are
complete for some intrinsic pseudo metric, as defined in \cite{FLW} and \cite{HKMW}.
\end{rem}
\begin{defi}\label{defi4.2}
A triangulation $\Tc=(\Kc,\Fc)$ is $\chi$-complete, if
\begin{itemize}
  \item[($C_1$)] $\Kc$ is $\chi$-complete.
  \item[($C_2$)] $\exists M>0,~\forall n\in \N,~e\in \Ec,$ such that
  $$\dfrac{1}{r(e)}\displaystyle\sum_{x\in\Fc_{e}}s(e,x)
  |d^{0}\chi_{n}(e^{-},x)+d^{0}\chi_{n}(e^{+},x)|^{2}\leq M.
$$
\end{itemize}
\end{defi}

For this type of 2-simplicial complexes one has
\begin{equation}
\forall p\in \N,~\exists n_{p},~n\geq n_{p}\Rightarrow \forall e\in \Ec,
\mbox{ such that }e^{+}\mbox{ or }e^{-}\in B_{p},~d^{0}\chi_{n}(e)=0.\label{prop1}
\end{equation}
\begin{equation}
\Ec=\displaystyle\bigcup_{n\in\N}\Ec_{n}\mbox{ if } \Ec_{n}:=\{e\in \Ec,~e^{+}\in B_{n}
\mbox{ or }e^{-}\in B_{n}\}.\label{prop2}
\end{equation}
\begin{equation}
\forall q\in \N,~\exists n_{q},~n\geq n_{q}\Rightarrow\forall(e,x)\in \Fc,
\mbox{ such that }e^-,~e^+\mbox{ or }x\in B_{q},~d^{0}\chi_{n}(e^{\pm},x)=0.\label{prop3}
\end{equation}
\begin{equation}
\Fc=\displaystyle\bigcup_{n\in\N}\Fc_{n}\mbox{ if } \Fc_{n}:=\{[x,y,z]\in \Fc,~x\in B_{n}
\mbox{ or }y\in B_{n}
\mbox{ or }z\in B_{n}\}.\label{prop4}
\end{equation}
\begin{equation}
\forall f\in l^2(\Vc),~\|f\|_{l^2(\Vc)}^2=\displaystyle\lim_{n\rightarrow\infty}
\langle\chi_{n}f,f\rangle_{l^2(\Vc)}.\label{prop5}
\end{equation}
\begin{equation}
\forall \varphi\in l^2(\Ec),~\|\varphi\|_{l^2(\Ec)}^2=\displaystyle\lim_{n\rightarrow\infty}
\frac{1}{2}\displaystyle\sum_{e\in\Ec}
r(e)\chi_{n}(e^{+})|\varphi(e)|^2.\label{prop6}
\end{equation}
\begin{equation}
\forall \phi\in l^2(\Fc),~\|\phi\|_{l^2(\Fc)}^2=\displaystyle\lim_{n\rightarrow\infty}
\frac{1}{6}\displaystyle\sum_{e\in\Ec}\widetilde{\chi_{n}}(e)\left(\displaystyle\sum_{x\in\Fc_{e}}
s(e,x)|\phi(e,x)|^2\right).\label{prop7}
\end{equation}
\begin{equation}
\displaystyle\lim_{n\rightarrow\infty}\displaystyle\sum_{e\in\Ec^*(n)}
r(e)|\varphi(e)|^2=0,\label{prop8}
\end{equation}
where
$$\Ec^*(n):=\{e\in\Ec,\exists x\in \Fc_e\mbox{ such that }
(e^{\pm},x)\in supp(d^0\chi_{n})\}
$$

\begin{equation}
\displaystyle\lim_{n\rightarrow\infty}
\displaystyle\sum_{e\in \Ec}\displaystyle\sum_{x\in \Fc^*_{e}(n)}s(e,x)|\phi(e,x)|^2=0,\label{prop9}
\end{equation}
where
$$\forall e\in\Ec,~\Fc^*_{e}(n):=\{x\in\Fc_e,~(e^{\pm},x)\in supp(d^{0}\chi_{n})\}.$$
\begin{propr}\label{youtube}
Let $\Tc$ be a simple triangulation of bounded degree, i.e
$\exists\lambda>0,~\forall x\in\Vc,~deg(x)\leq\lambda.$
Then $\Tc$ is a $\chi$-complete triangulation.
\end{propr}
\begin{dem}

Let us consider $\Tc$ an infinite triangulation. Given $o\in \Vc,$ let $B_{n}$
be a ball of radius $n\in \N$ centered by the vertex $o$:
$$B_{n}=\{x\in \Vc,~d_{comb}(o,x)\leq n\}.
$$

We set the cut-off function $\chi_{n}\in \Cc_{c}(\Vc)$ as follow:
$$\chi_{n}(x):=\left(\dfrac{2n-d_{comb}(o,x)}{n}\vee 0 \right)\wedge 1,~\forall n\in \N^{*}.
$$
\begin{enumerate}
  \item[-] If $x\in B_{n}\Rightarrow \chi_{n}(x)=1$ and $x\in B_{2n}^{c}\Rightarrow \chi_{n}(x)=0.$
  \item[-] For $e\in \Ec,$ we have that
$$|d^0\chi_{n}(e)|\leq\dfrac{1}{n}\left|d_{comb}(o,e^{+})-d_{comb}(o,e^{-})\right|=\dfrac{1}{n}.
$$

Hence
$$\forall x\in \Vc,~\displaystyle\sum_{e\in\Ec,e^{\pm}=x}
|d^{0}\chi_{n}(e)|^{2}\leq \dfrac{\lambda}{n^2}
$$

and
$$\forall e\in \Ec,~\displaystyle\sum_{x\in\Fc_e}
|d^{0}\chi_{n}(e^{-},x)+d^{0}\chi_{n}(e^{+},x)|^{2}
\leq\dfrac{2\lambda}{n^2}.
$$
\end{enumerate}
\end{dem}
\begin{exa}\emph{(A $\chi$-complete triangulation)}

We consider $\Tc$ a 6-regular simple triangulation, i.e. $deg(x)=6,~\forall x\in\Vc.$
Then, by Proposition \ref{youtube} we have that $\Tc$ is a $\chi$-complete triangulation.
\begin{figure}[ht]
\includegraphics*[width=7cm]{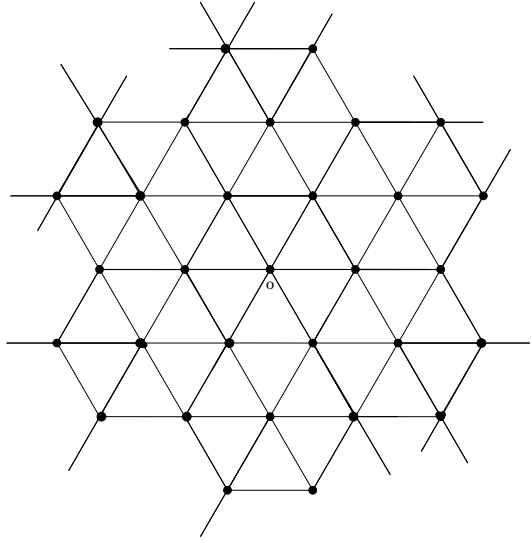}
\caption{\label{Paving} An infinite 6-regular triangulation}
\end{figure}
\end{exa}

\begin{propr}
Let $\Tc=(\Kc$,$\Fc)$ be a $\chi$-complete triangulation. Then
$$Dom\left((L_0)_{min}\right)=Dom(\delta_{min}^0d_{min}^0).
$$
\end{propr}
\begin{dem}

In Proposition \ref{Yasch}, we have already $(L_0)_{min}\subseteq \delta_{min}^0d_{min}^0.$
Indeed, we will show that $\delta_{min}^0d_{min}^0\subseteq(L_0)_{min}.$ Let $f\in Dom(\delta_{min}^0d_{min}^0),$
by the $\chi$-completeness of $\Tc,$ we now consider a sequence $(\chi_nf)_n\subseteq\Cc_{c}(\Vc).$
It remains to show that:
\begin{equation}\label{en}
\displaystyle\lim_{n\rightarrow\infty}\|f-\chi_{n}f\|_{l^2(\Vc)}
+\|L_0(f-\chi_{n}f)\|_{l^2(\Vc)}=0.
\end{equation}

For the first term of (\ref{en}), since $f\in l^2(\Vc)$ we have
$$\|f-\chi_{n}f\|^2_{l^2(\Vc)}\leq\displaystyle\sum_{x\in B^c_n}c(x)|f(x)|^2\rightarrow 0,
\mbox{ when }n\rightarrow\infty.
$$

For the second term of (\ref{en}), we need a derivation formula of $d^0,$ see \cite{M2}.
Let $e\in\Ec,$ for each
$(f,g)\in\Cc_{c}(\Vc)\times\Cc_{c}(\Vc)$ we have
\begin{equation}\label{eqqq}
d^0(fg)(e)=f(e^+)d^0(g)(e)+d^0(f)(e)g(e^-).
\end{equation}

By the definition of $L_0,$ we have
$$\|L_0(f-\chi_{n}f)\|^2_{l^2(\Vc)}=\displaystyle\sum_{x\in\Vc}
\dfrac{1}{c(x)}\left|\displaystyle\sum_{e,e^+=x}r(e)d^0((1-\chi_n)f)(e)\right|^2.
$$

Using the derivation formula (\ref{eqqq}), we get
\begin{equation*}
\begin{split}
\|L_0(f-\chi_{n}f)\|^2_{l^2(\Vc)}
&\leq2\displaystyle\sum_{x\in\Vc}\dfrac{1}{c(x)}\left|\displaystyle\sum_{e,e^+=x}
r(e)(1-\chi_n)(e^+)d^0(f)(e)\right|^2\\
&+2\displaystyle\sum_{x\in\Vc}\dfrac{1}{c(x)}
\left|\displaystyle\sum_{e,e^+=x}r(e)f(e^-)d^0(\chi_n)(e)\right|^2\\
&=2\left(\|(1-\chi_n)L_0(f)\|^2_{l^2(\Vc)}+\displaystyle\sum_{x\in\Vc}\dfrac{1}{c(x)}
\left|\displaystyle\sum_{e,e^+=x}r(e)f(e^-)d^0(\chi_n)(e)\right|^2\right).
\end{split}
\end{equation*}

Since $L_0(f)\in l^2(\Vc)$, we have
$$\displaystyle\lim_{n\rightarrow\infty}\|(1-\chi_n)L_0(f)\|_{l^2(\Vc)}=0.
$$

On the other hand, by the hypothesis iii) of $\chi$-completeness and the Cauchy-Schwarz
inequality, we get
\begin{equation*}
\begin{split}
\displaystyle\sum_{x\in\Vc}\dfrac{1}{c(x)}
\left|\displaystyle\sum_{e,e^+=x}r(e)f(e^-)d^0(\chi_n)(e)\right|^2
&\leq\displaystyle\sum_{x\in\Vc}\dfrac{1}{c(x)}
\left(\displaystyle\sum_{e,e^+=x}r(e)|d^0(\chi_n)(e)|^2\right)\\
&\left(\displaystyle\sum_{e\in supp(d^0\chi_n),e^+=x}r(e)|f(e^-)|^2\right)\\
&\leq\displaystyle\sum_{x\in\Vc}C\displaystyle\sum_{e\in supp(d^0\chi_n),e^+=x}r(e)|f(e^-)|^2\\
&\leq C\displaystyle\sum_{e\in supp(d^0\chi_n)}r(e)|f(e^-)|^2.
\end{split}
\end{equation*}

The properties (\ref{prop1}) and (\ref{prop2}) permit to conclude that this term
tends to $0$ when $\infty.$
\end{dem}

\subsection{The case of a not $\chi$-complete triangulation }

In \cite{BGJ}, the authors use the offspring function on the trees to give
a counter example of a graph which is not $\chi$-complete. The same thing for the
triangulations is not always $\chi$-complete. To prove it, we will study
the triangular tree in Definition \ref{def.}.

Let $\Tc$ a weighted triangulation, one can take any point $o\in\Vc.$
Given $n\in\N,$ we denote the spheres by
$$\Sc_n:=\{x\in \Vc,~d_{comb}(o,x)=n\}.$$

\begin{defi}\label{def.}
A triangular tree $\Tc=(\Vc,\Fc)$ with the origin vertex $o$ is a triangulation where
$\Vc=\displaystyle\cup_{n\in\N}S_n,$ such that
$$\forall x\in\Sc_n\backslash\{o\},~\Vc(x)\cap\Sc_{n-1}=\{\overleftarrow{x}\}.
$$
$$\forall x\in\Sc_n,~y\in\Vc(x)\cap\Sc_{n+1}\Leftrightarrow\overleftarrow{y}=x.
$$
$$(x,y)\in\Ec\cap(\Sc_n\backslash\{o\})^2\Rightarrow\overleftarrow{x}=\overleftarrow{y}.
$$
where $\overleftarrow{x}$ the unique vertex in $\Sc_{n-1},$ which is related with
$x\in\Sc_n\backslash\{o\}.$

\begin{figure}[ht]
\includegraphics*[width=7cm]{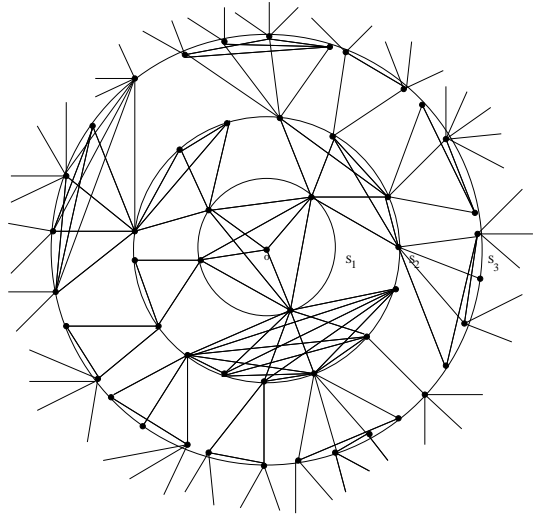}
\caption{\label{Tree} A Triangular Tree}
\end{figure}
\end{defi}

Let $\Tc$ be a simple triangular tree. \emph{The offspring} of the \emph{n}-th generation (see \cite{BGJ}) is given by
$$\mbox{off}(n)=\dfrac{\#\Sc_{n+1}}{\#\Sc_n}.
$$
\begin{propr}
Let $\Tc$ be a simple triangular tree with the origin vertex $o.$ Assume that
$$\displaystyle\sup_{n\in\N}\displaystyle\sup_{x\in\Sc_n}
\dfrac{\#\left(\Vc(x)\cap\Sc_{n+1}\right)}{\emph{off}(n)}<\infty.
$$

Then
$$\Tc\mbox{ is }\chi\mbox{-complete}\Leftrightarrow\displaystyle\sum_{n\geqslant 1}
\dfrac{1}{\sqrt{\emph{off}(n)}}=\infty.
$$
\end{propr}
\begin{dem}
\begin{enumerate}
\item[$\Rightarrow)$] In a proof by contradiction, we start by assuming that
$\Tc$ is $\chi$-complete and the series converges. So, there exists a sequence
$(\chi_n)_n$ included in $\Cc_{c}(\Vc),$ satisfying
$$\exists C>0,~\forall n\in\N,~\displaystyle\sum_{y\sim x}
|\chi_n(x)-\chi_n(y)|^2\leq C,~x\in\Vc.
$$

Given $n,m\in\N$ and $x_m\in\Sc_{m}.$ By the local finiteness of the triangulation,
we find $x_{m+1}\in\Vc(x_m)\cap\Sc_{m+1},$ such that
$$|\chi_{n}(x_m)-\chi_{n}(x_{m+1})|=\displaystyle\min_{y\in\Vc(x_m)\cap\Sc_{m+1}}|\chi_{n}(x_m)-\chi_{n}(y)|.
$$

But,
$$\displaystyle\sum_{y\in\Vc(x_m)\cap\Sc_{m+1}}|\chi_{n}(x_m)-\chi_{n}(y)|^2\leq C.
$$

Hence
$$|\chi_{n}(x_m)-\chi_{n}(x_{m+1})|\leq\dfrac{\sqrt{C}}{\sqrt{\mbox{off}(m)}}.
$$

Moreover, by convergence of the series, there is $N\in\N$ such that
$$\displaystyle\sum_{k\geq N}\dfrac{1}{\sqrt{\mbox{off}(k)}}
<\dfrac{1}{2\sqrt{C}}.
$$

Then, by ii) of the definition of $\chi$-completeness, there is $n_0\in\N$ such that
$\chi_{no}(x)=1$ for all $d_{comb}(o,x)\leq N.$ Since $\chi_{no}$
is with finite support, there is $M\in\N$ such that $\chi_{no}(x)=0$
for all $d_{comb}(o,x)\geq N+M.$ Therefore,
\begin{equation*}
\begin{split}
|\chi_{no}(x_N)-\chi_{no}(x_{N+M})|
& \leq|\chi_{no}(x_N)-\chi_{no}(x_{N+1})|+.....
+|\chi_{no}(x_{N+M-1})-\chi_{no}(x_{N+M})|\\
& \leq\sqrt{C}\displaystyle\sum_{k=n}^{N+M-1}
\dfrac{1}{\sqrt{\mbox{off}(k)}}<\frac{1}{2}.
\end{split}
\end{equation*}
Since $|\chi_{no}(x_N)-\chi_{no}(x_{N+M})|=1,$ we have the contradiction.
\item[$\Leftarrow)$]We consider the cut-off function:
$$ \chi_n(x)=\left\{\begin{array}{ll}
                        1 & \hbox{if  }~ d_{comb}(o,x)\leqslant n, \\
                        \max\left(0,1-\displaystyle\sum_{k=n}^{d_{comb}(o,x)-1}
                        \dfrac{1}{\sqrt{\mbox{off}(k)}}\right)
                        & \hbox{if  } \, d_{comb}(o,x)> n.

                 \end{array}
            \right.
   $$

Since the series diverges, $\chi_n$ is with finite support and satisfies
i) and ii) of the definition of $\chi$-completeness. Given $x\in\Sc_m$ with $m>n,$ we have
\begin{equation*}
\displaystyle\sum_{y\in\Vc(x)\cap\Sc_{m+1}}|\chi_n(x)-\chi_n(y)|^2
\leq\dfrac{\#\left(\Vc(x)\cap\Sc_{m+1}\right)}{\mbox{off}(m)}.
\end{equation*}
\begin{equation*}
\displaystyle\sum_{y\in\Vc(x)\cap\Sc_m}|\chi_n(x)-\chi_n(y)|^2=0.
\end{equation*}
\begin{equation*}
\displaystyle\sum_{y\in\Vc(x)\cap\Sc_{m-1}}|\chi_n(x)-\chi_n(y)|^2
=|\chi_n(x)-\chi_n(\overleftarrow{x})|^2
\leq\dfrac{1}{\mbox{off}(m-1)}.
\end{equation*}

On the other hand,
\begin{enumerate}
\item[i)]If $e\in \Sc_m\times\Sc_{m+1}$ with $m>n,$ we have
\begin{equation*}
\begin{split}
\displaystyle\sum_{x\in \Fc_e}
|d^{0}\chi_{n}(e^{-},x)+d^{0}\chi_{n}(e^{+},x)|^{2}
& =\displaystyle\sum_{x\in \Fc_e}|2\chi_{n}(x)-\chi_{n}(e^{-})-\chi_{n}(e^{+})|^{2}\\
& \leq\dfrac{|\Fc_{e}|}{\mbox{off}(m)}\leq\dfrac{\#\left(\Vc(e^-)\cap\Sc_{m+1}\right)}{\mbox{off}(m)}.
\end{split}
\end{equation*}
\item[ii)]If $e\in \Sc_m\times\Sc_m$ with $m>n,$ we have
\begin{equation*}
\begin{split}
\displaystyle\sum_{x\in \Fc_e}
|d^{0}\chi_{n}(e^{-},x)+d^{0}\chi_{n}(e^{+},x)|^{2}
& =\displaystyle\sum_{x\in \Fc_e}
|2\chi_{n}(x)-\chi_{n}(e^{-})-\chi_{n}(e^{+})|^{2}\\
& =|2\chi_{n}(\overleftarrow{e})-\chi_{n}(e^{-})-\chi_{n}(e^{+})|^{2}\\
& \leq \dfrac{4}{\mbox{off}(m-1)},
\end{split}
\end{equation*}
with $\overleftarrow{e}$ is a unique vertex in $\Sc_{m-1}\cap\Fc_e.$
\end{enumerate}
It satisfies Definition \ref{defi4.2} of $\chi$-completeness.
\end{enumerate}
\end{dem}

\begin{coro}
Let $\Tc$ be a simple triangular tree, endowed with an origin such that
$$\emph{off}(n)=\#\left(\Vc(x)\cap\Sc_{n+1}\right),\mbox{ for all } x\in \Sc_n,
$$
then $\Tc$ is $\chi$-complete if and only if
$$\displaystyle\sum_{n\geqslant 1}\dfrac{1}{\sqrt{\emph{off}(n)}}=\infty.
$$
\end{coro}

\begin{exa}
Set $\alpha>0.$ Let $\Tc$ be a simple triangular tree, endowed with an origin such that
$$\emph{off}(n)=\#\left(\Vc(x)\cap\Sc_{n+1}\right)=\lfloor n^\alpha \rfloor+1,\mbox{ for all } x\in \Sc_n,
$$
then $\Tc$ is $\chi$-complete if only if $\alpha\leq2.$
\end{exa}

\section{Essential self-adjointness}
In \cite{AT}, the authors use the $\chi$-completeness hypothesis on a graph
to ensure essential self-adjointness for the Gau$\beta$-Bonnet operator and the Laplacian.
In this section, with the same idea we will prove the main result, when the triangulation
is $\chi$-complete. Let us begin from
\begin{propr} \label{pro2}
Let $\Tc=(\Kc,\Fc)$ be a $\chi$-complete triangulation then the operator $d^1+\delta^1$
is essentially self-adjoint on $\Cc_{c}(\Ec)\oplus\Cc_{c}(\Fc).$
\end{propr}
\begin{dem}
It suffices to show that $d^1_{min}=d^1_{max}$ and $\delta^1_{min}=\delta^1_{max}.$
Indeed, $d^1+\delta^1$ is a direct sum and if $F=(\varphi,\phi)\in Dom((d^1+\delta^1)_{max})$
then $\varphi\in Dom(d^1_{max})$ and $\phi\in Dom(\delta^1_{max}).$ By hypothesis, we have
$\varphi\in Dom(d^1_{min})$ and $\phi\in Dom(\delta^1_{min}),$ thus $F\in Dom((d^1+\delta^1)_{min}).$
\begin{enumerate}
\item[1)]Let $\varphi\in Dom(d^1_{max}),$ we will show that
$$\|\varphi-\widetilde{\chi_{n}}\varphi\|_{l^2(\Ec)}
+\|d^1\left(\varphi-\widetilde{\chi_{n}}\varphi\right)\|_{l^2(\Fc)}
\rightarrow 0\mbox{ when }n\rightarrow\infty.
$$

By the properties (\ref{prop1}) and (\ref{prop2}), we know that
$$\forall p\in \N,~\exists n_p,~\forall n\geq n_p,~
\| \varphi-\widetilde{\chi_{n}}\varphi\|^2_{l^2(\Ec)}
\leq\displaystyle\sum_{e\in\Ec_p^{c}}r(e)|\varphi(e)|^2
$$
so $\displaystyle\lim_{n\longrightarrow\infty}
\|\varphi-\widetilde{\chi_{n}}\varphi\|=0.$

From the derivation formula (\ref{Eq4.1}) in Lemma \ref{lem3}, we have
\begin{equation*}
\begin{split}
d^1\left(\varphi-\widetilde{\chi_{n}}\varphi\right)(e,x)
& = d^1\left(\left(\widetilde{1-\chi_{n}}\right)\varphi\right)(e,x)\\
& = \left(1-\widetilde{\widetilde{\chi_{n}}}\right)(e,x)
d^1(\varphi)(e,x)\\
& + \frac{1}{6}\left(d^{0}(1-\chi_{n})(x,e^-)+d^{0}(1-\chi_{n})(x,e^+)\right)\varphi(e)\\
& + \frac{1}{6}\left(d^{0}(1-\chi_{n})(e)+d^{0}(1-\chi_{n})(e^-,x)\right)\varphi(e^+,x)\\
& + \frac{1}{6}\left(d^{0}(1-\chi_{n})(e^+,x)+d^{0}(1-\chi_{n})(-e)\right)\varphi(x,e^-)\\
& = \left(1-\widetilde{\widetilde{\chi_{n}}}\right)(e,x)
d^1(\varphi)(e,x)\\
& + \frac{1}{6}\left(d^{0}\chi_{n}(e^-,x)+d^{0}\chi_{n}(e^+,x)\right)\varphi(e)\\
& + \frac{1}{6}\left(d^{0}\chi_{n}(-e)+d^{0}\chi_{n}(x,e^-)\right)\varphi(e^+,x)\\
& + \frac{1}{6}\left(d^{0}\chi_{n}(x,e^+)+d^{0}\chi_{n}(e)\right)\varphi(x,e^-).
\end{split}
\end{equation*}
Since $d^1\varphi\in l^2(\Fc),$ one has
$$\displaystyle\lim_{n\rightarrow\infty}
\|\left(1-\widetilde{\widetilde{\chi_{n}}}\right)d^1\varphi\|_{l^2(\Fc)}=0.
$$
On the other hand,
\begin{equation*}
\begin{split}
\displaystyle\sum_{(e,x)\in\Fc}s(e,x)|\varphi(e)|^2
|d^{0}\chi_{n}(e^-,x)+d^{0}\chi_{n}(e^+,x)|^2
& = \displaystyle\sum_{e\in\Ec}|\varphi(e)|^2\displaystyle\sum_{x\in\Fc_e}
s(e,x)|d^{0}\chi_{n}(e^-,x)+d^{0}\chi_{n}(e^+,x)|^2\\
& \leq M \displaystyle\sum_{e\in\Ec^*(n)}r(e)|\varphi(e)|^2.
\end{split}
\end{equation*}
The property (\ref{prop8}) allows to conclude that this term tends to $0$ as $n\rightarrow\infty.$
Applying the same process to the other terms, we have
\begin{equation*}
\begin{split}
\displaystyle\sum_{(e,x)\in\Fc}s(e,x)|\varphi(e^+,x)|^2
|d^{0}\chi_{n}(-e)+d^{0}\chi_{n}(x,e^-)|^2
& = \displaystyle\sum_{(e^+,x)\in\Ec}|\varphi(e^+,x)|^2
\displaystyle\sum_{y\in\Fc_{(e^+,x)}}s(e^+,x,y)\\
& |d^{0}\chi_{n}(e^+,y)+d^{0}\chi_{n}(x,y)|^2\\
& \leq M \displaystyle\sum_{(e^+,x)\in\Ec^*(n)}r(e^+,x)|\varphi(e^+,x)|^2
\end{split}
\end{equation*}
and
\begin{equation*}
\begin{split}
\displaystyle\sum_{(e,x)\in\Fc}s(e,x)|\varphi(x,e^-)|^2
|d^{0}\chi_{n}(x,e^+)+d^{0}\chi_{n}(e)|^2
& = \displaystyle\sum_{(x,e^-)\in\Ec}|\varphi(x,e^-)|^2
\displaystyle\sum_{y\in\Fc_{(x,e^-)}}s(x,e^-,y)\\
& |d^{0}\chi_{n}(x,y)+d^{0}\chi_{n}(e^-,y)|^2\\
& \leq M \displaystyle\sum_{(x,e^-)\in\Ec^*(n)}r(x,e^-)|\varphi(x,e^-)|^2.
\end{split}
\end{equation*}
\item[2)] Let $\phi\in Dom(\delta^1_{max}),$ we will show that
$$\|\phi-\widetilde{\widetilde{\chi_{n}}}\phi\|_{l^2(\Fc)}
+\|\delta^1(\phi-\widetilde{\widetilde{\chi_{n}}}\phi)\|_{l^2(\Ec)}
\rightarrow 0\mbox{ when }n\rightarrow\infty.
$$
By the properties (\ref{prop3}) and (\ref{prop4}), we know that
\begin{equation*}
\begin{split}
\|\phi-\widetilde{\widetilde{\chi_{n}}}\phi\|^2_{l^2(\Fc)}
& = \frac{1}{6}\displaystyle\sum_{(x,y,z)\in\Fc}s(x,y,z)|1-\widetilde{\widetilde{\chi_{n}}}(x,y,z)|^2
|\phi(x,y,z)|^2\\
& \leq  \displaystyle\sum_{(x,y,z)\in\Fc^c_q}s(x,y,z)|\phi(x,y,z)|^2
\rightarrow 0,\mbox{ when } n\rightarrow\infty.
\end{split}
\end{equation*}
By the derivation formula (\ref{Eq4.2}) in Lemma \ref{lem3}, we have
\begin{equation*}
\begin{split}
\delta^1(\phi-\widetilde{\widetilde{\chi_{n}}}\phi)(e)
& = \delta^1\left((\widetilde{\widetilde{1-\chi_{n}}})\phi\right)(e)\\
& = (1-\widetilde{\chi_{n}})(e)\delta^1(\phi)(e)\\
& + \dfrac{1}{6r(e)}\displaystyle\sum_{x\in\Fc_e}s(e,x)
d^0(1-\chi_{n})(e^-,x)\phi(e,x)\\
& + \dfrac{1}{6r(e)}\displaystyle\sum_{x\in\Fc_e}s(e,x)
d^0(1-\chi_{n})(e^+,x)\phi(e,x)\\
& = \left(1-\widetilde{\chi_{n}}\right)(e)\delta^1(\phi)(e)\\
& + \dfrac{1}{6r(e)}\displaystyle\sum_{x\in\Fc_e}s(e,x)
d^0(\chi_{n})(x,e^-)\phi(e,x)\\
& + \dfrac{1}{6r(e)}\displaystyle\sum_{x\in\Fc_e}s(e,x)
d^0(\chi_{n})(x,e^+)\phi(e,x).
\end{split}
\end{equation*}
We know that
$$\displaystyle\lim_{n\rightarrow\infty}\|\left(1-\widetilde{\chi_{n}}\right)
\delta^1(\phi)\|=0
$$
because $\delta^1\phi\in l^2(\Ec).$ For the second and third terms,
we use the inequality of Definition \ref{defi4.2} and the Cauchy-Schwarz inequality.
Fix $e\in \Ec,$ then
\begin{equation*}
\begin{split}
\left|\displaystyle\sum_{x\in \Fc_{e}}s(e,x)
\left(d^0(\chi_{n})(x,e^-)+d^0(\chi_{n})(x,e^+)\right)
\phi(e,x)\right|^2
& \leq \displaystyle\sum_{x\in \Fc_{e}}s(e,x)|d^0(\chi_{n})(x,e^-)+d^0(\chi_{n})(x,e^+)|^2\\
& \times \displaystyle\sum_{x\in\Fc_{e}^*(n)}s(e,x)|\phi(e,x)|^2\\
& \leq Mr(e)\displaystyle\sum_{x\in\Fc_{e}^*(n)}s(e,x)|\phi(e,x)|^2.
\end{split}
\end{equation*}
Therefore,
\begin{equation*}
\begin{split}
\displaystyle\sum_{e\in \Ec}r(e)\left|\dfrac{1}{r(e)}
\displaystyle\sum_{x\in \Fc_{e}}s(e,x)
\left(d^0(\chi_{n})(x,e^-)+d^0(\chi_{n})(x,e^+)\right)\phi(e,x)\right|^2
& \leq M\displaystyle\sum_{e\in \Ec}\displaystyle
\sum_{x\in\Fc_{e}^*(n)}s(e,x)|\phi(e,x)|^2.
\end{split}
\end{equation*}
By property (\ref{prop9}), this terms tends to $0.$
\end{enumerate}
\end{dem}
\begin{coro}
Let $\Tc=(\Kc,\Fc)$ be a $\chi$-complete triangulation then the operator $L_1^{+}\oplus L_2$
is essentially self-adjoint on $\Cc_{c}(\Ec)\oplus\Cc_{c}(\Fc).$
\end{coro}
\begin{dem}
First we have that $L_1^{+}\oplus L_2=\left(d^1+\delta^1\right)^2$ and
$L_1^{+}\oplus L_2\left(\Cc_{c}(\Ec)\oplus\Cc_{c}(\Fc)\right)\subseteq
\Cc_{c}(\Ec)\oplus\Cc_{c}(\Fc).$ As Proposition 13 in \cite{AT} we
prove that $d^1+\delta^1$ is essentially self-adjoint if and only if
$L_1^{+}\oplus L_2$ is essentially self-adjoint.
\end{dem}
\begin{thm}\label{esthm}
Let $\Tc=(\Kc,\Fc)$ be a $\chi$-complete triangulation then the operator
$T$ is essentially self-adjoint on $\Cc_{c}(\Vc)\oplus\Cc_{c}(\Ec)\oplus\Cc_{c}(\Fc).$
\end{thm}
\begin{dem}
\emph{\underline{$1^{st}$ Step:}} We will show that
$$Dom(T_{min})=Dom(d^{0}_{min})\oplus\left(Dom(\delta^{0}_{min})\cap Dom(d^{1}_{min})\right)
\oplus Dom(\delta^{1}_{min}).
$$

Let $F=(f,\varphi,\phi)\in Dom(d^{0}_{min})\oplus
\left(Dom(\delta^{0}_{min})\cap Dom(d^{1}_{min})\right)\oplus Dom(\delta^{1}_{min}).$
Then there exist $(f_n)_n\subseteq\Cc_{c}(\Vc)$ and
$(\phi_n)_n\subseteq\Cc_{c}(\Fc)$ such that:
\begin{enumerate}
\item[-]$f_n\rightarrow f\mbox{ in }l^2(\Vc)$ and
$d^0 f_n \rightarrow d^0_{min} f\mbox{ in }l^2(\Ec).$
\item[-]$\phi_n\rightarrow\phi\mbox{ in }l^2(\Fc)$ and
$\delta^1 \phi_n\rightarrow\delta^1_{min}\phi\mbox{ in }l^2(\Ec)$
\end{enumerate}

On the other hand, let $\varphi\in Dom(\delta^{0}_{min})\cap Dom(d^{1}_{min}).$
By the $\chi$-completeness of $\Tc,$ we now consider the sequence
$(\widetilde{\chi_{n}}\varphi)_n\subseteq\Cc_{c}(\Ec).$
It remains to show that
$$\|\varphi-\widetilde{\chi_{n}}\varphi\|_{l^2(\Ec)}
+\|\delta^0(\varphi-\widetilde{\chi_{n}}\varphi)\|_{l^2(\Vc)}
+\|d^1(\varphi-\widetilde{\chi_{n}}\varphi)\|_{l^2(\Fc)}\rightarrow 0,
\mbox{ when } n\rightarrow\infty.
$$

The first and the third terms has already been shown in Proposition \ref{pro2}.
For the following we need a derivation formula of $\delta^0$ taken in \cite{M2}.
Let $x\in\Vc,$ for each $(f,\varphi)\in\Cc_{c}(\Vc)\times\Cc_{c}(\Ec)$
we have
\begin{equation}\label{eqq}
\delta^0(\widetilde{f}\varphi)(x)=f(x)\delta^0(\varphi)(x)
-\dfrac{1}{2c(x)}\displaystyle\sum_{e,e^+=x}r(e)d^0(f)(e)\varphi(e).
\end{equation}

Therefore, by derivation formula (\ref{eqq}), we get
$$\delta^0(\varphi-\widetilde{\chi_{n}}\varphi)(x)
=(1-\chi_{n})(x)\delta^0(\varphi)(x)
+\dfrac{1}{2c(x)}\displaystyle\sum_{e,e^+=x}r(e)d^0\chi_{n}(e)\varphi(e).
$$

As a consequence, because $\delta^0\varphi\in l^2(\Vc),$ we have
$$\lim\limits_{n\rightarrow\infty}\|(1-\chi_{n})\delta^0\varphi\|_{l^2(\Vc)}=0.
$$

For the second term, we combine the property iii) of $\chi$-completeness for a graph
with the Cauchy-Schwarz inequality to obtain for all $x\in\Vc,$
\begin{equation*}
\begin{split}
|\displaystyle\sum_{e,e^+=x}r(e)d^0\chi_n(e)\varphi(e)|^2
&\leq\displaystyle\sum_{e,e^+=x}r(e)|d^0\chi_n(e)|^2\displaystyle\sum_{e\in supp(d^0\chi_n),e^+=x}
r(e)|\varphi(e)|^2\\
&\leq Cc(x)\displaystyle\sum_{e\in supp(d^0\chi_n),e^+=x}r(e)|\varphi(e)|^2.
\end{split}
\end{equation*}

So,
\begin{equation*}
\begin{split}
\displaystyle\sum_{x\in\Vc}\dfrac{1}{c(x)}|\displaystyle\sum_{e,e^+=x}r(e)d^0\chi_n(e)\varphi(e)|^2
&\leq\displaystyle\sum_{x\in\Vc}C\displaystyle\sum_{e\in supp(d^0\chi_n),e^+=x}
r(e)\varphi(e)|^2\\
&\leq C\displaystyle\sum_{e\in supp(d^0\chi_n)}
r(e)|\varphi(e)|^2\rightarrow 0,~\mbox{ when }n\rightarrow\infty,
\end{split}
\end{equation*}

by the properties (\ref{prop1}) and (\ref{prop2}).

Hence
$$F_n\rightarrow F\mbox{ in }\Hc,~TF_n\rightarrow T_{min}F
\mbox{ in }\Hc,
$$
where $F_n=(f_n,\widetilde{\chi}_n\varphi,\phi_n)
\mbox{ and }T_{min}F(f,\varphi,\phi)=(\delta^0_{min} \varphi,d^0_{min}f
+\delta^1_{min} \phi,d^1_{min} \varphi).$ Then $F\in Dom(T_{min}).$

\emph{\underline{$2^{th}$ Step:}} To show that $T$ is essentially self-adjoint,
we will prove that $T_{max}=T_{min}.$ By the first step, Theorem 1 in \cite{AT}
and Proposition \ref{pro2} it remains to show that:
$$Dom(T_{max})\subseteq Dom(d^{0}_{max})
\oplus\left(Dom(\delta^{0}_{max})\cap Dom(d^{1}_{max})\right)
\oplus Dom(\delta^{1}_{max}).
$$

Let $F=(f,\varphi,\phi)\in Dom(T_{max})$ then $TF\in\Hc.$ This implies that
$\delta^0\varphi\in l^2(\Vc),$ $d^0f+\delta^1\phi\in l^2(\Ec)$ and $d^1\varphi\in l^2(\Fc).$
As consequence, by the definition of $\delta^{0}_{max}$ and $d^{1}_{max}$
we have $\varphi\in Dom(\delta^{0}_{max})\cap Dom(d^{1}_{max}).$ Moerever, by $\chi$-completness
of $\Tc,$ there exists a sequence of cut-off functions $\left(\chi_n\right)_n\subseteq \Cc_{c}(\Vc).$
Then, the parallelogram identity with Lemma \ref{lem2} we get
$$\|d^0(\chi_n f)+\delta^1(\widetilde{\widetilde{\chi_n}}\phi)\|^2_{l^2(\Ec)}
=\|d^0\chi_nf\|^2_{l^2(\Ec)}+\|\delta^1\widetilde{\widetilde{\chi_n}}\phi\|^2_{l^2(\Ec)}.
$$

Now, it remains to prove that $d^0(\chi_n f)+\delta^1(\widetilde{\widetilde{\chi_n}}\phi)$
converges in $l^2(\Ec).$ Indeed, we need some formulas taken in Lemma \ref{lem3} and \cite{M2}
to give that:
$$d^0(\chi_nf)=\widetilde{\chi_n}d^0(f)+\widetilde{f}d^0(\chi_n).
$$
$$\delta^1(\widetilde{\widetilde{\chi_n}}\phi)(e)=\widetilde{\chi_n}(e)\delta^1(\phi)(e)
+\underbrace{\dfrac{1}{6r(e)}\displaystyle\sum_{x\in\Fc_e}s(e,x)
\left[d^0(\chi_n)(e^-,x)+d^0(\chi_n)(e^+,x)\right]\phi(e,x)}_{\Ic_n(e)}.
$$

Therefore, we have
\begin{equation*}
\begin{split}
\|d^0(f-\chi_nf)+\delta^1(\phi-\widetilde{\widetilde{\chi_n}}\phi)\|^2_{l^2(\Ec)}
&=\|(1-\widetilde{\chi_n})(d^0f+\delta^1\phi)+\widetilde{f}d^0\chi_n+\Ic_n\|^2_{l^2(\Ec)}\\
&\leq3\left(\|(1-\widetilde{\chi_n})(d^0f+\delta^1\phi)\|^2_{l^2(\Ec)}
+\|\widetilde{f}d^0\chi_n\|^2_{l^2(\Ec)}+\|\Ic_n\|^2_{l^2(\Ec)}\right)
\end{split}
\end{equation*}

Because $d^0f+\delta^1\phi\in l^2(\Ec),$ we have
$$\displaystyle\lim_{n\rightarrow\infty}
\|(1-\widetilde{\chi_n})(d^0f+\delta^1\phi)\|^2_{l^2(\Ec)}=0.
$$

By Proposition \ref{pro2} we have
$$\displaystyle\lim_{n\rightarrow\infty}\|\Ic_n\|^2_{l^2(\Ec)}=0.
$$

Moreover, by the hypothesis iii) of $\chi$-completeness we have
\begin{equation*}
\begin{split}
\|\widetilde{f}d^0\chi_n\|^2_{l^2(\Ec)}
&=\dfrac{1}{2}\displaystyle\sum_{e\in\Ec}r(e)|\widetilde{f}(e)d^0(\chi_n)(e)|^2\\
&\leq\displaystyle\sum_{e\in\Ec}r(e)|f(e^+)|^2|d^0(\chi_n)(e)|^2\\
&=\displaystyle\sum_{x\in\Vc}|f(x)|^2\displaystyle\sum_{e,e^+=x}|r(e)d^0\chi_n(e)|^2\\
&\leq C\displaystyle\sum_{x\in\Vc_n}c(x)|f(x)|^2
\end{split}
\end{equation*}
where $\Vc_n:=\{x\in\Vc,\exists e\in supp(d^0\chi_n)\mbox{ such that }e^+=x\}.$
This term tends to $0$ by the property (\ref{prop2}).
\end{dem}
\begin{thm}
Let $\Tc=(\Kc,\Fc)$ be a $\chi$-complete triangulation. Then
$T$ is essentially self-adjoint on $\Cc_{c}(\Vc)\oplus\Cc_{c}(\Ec)\oplus\Cc_{c}(\Fc)$
if and only if
$L$ is essentially self-adjoint on $\Cc_{c}(\Vc)\oplus\Cc_{c}(\Ec)\oplus\Cc_{c}(\Fc).$
\end{thm}
\begin{dem}

Since
$$T(\Cc_{c}(\Vc)\oplus\Cc_{c}(\Ec)\oplus\Cc_{c}(\Fc))
\subseteq\Cc_{c}(\Vc)\oplus\Cc_{c}(\Ec)\oplus\Cc_{c}(\Fc),
$$

using the same technique in the proof of Proposition 13 in \cite{AT}, the result holds.
\end{dem}
\begin{coro}\label{co5}
Let $\Tc=(\Kc,\Fc)$ be a $\chi$-complete triangulation then
$L$ is essentially self-adjoint on $\Cc_{c}(\Vc)\oplus\Cc_{c}(\Ec)\oplus\Cc_{c}(\Fc).$
\end{coro}

\section{Examples}
\subsection{A triangulation with 1-dimensional decomposition}

We now strengthen the previous example and follow ideas of \cite{BG} and \cite{BGJ}.
\begin{defi}\emph{(1-dimensional decomposition)} A \emph{1-dimensional decomposition}
of the graph $\Kc=(\Vc,\Ec)$ is a family of finite sets $(\Sc_n)_{n\in\N}$ which forms
a partition of $\Vc,$ that is $\Vc=\displaystyle\sqcup_{n\in\N}\Sc_n,$ such that for all
$x\in\Sc_n,y\in\Sc_m,$
$$(x,y)\in\Ec\Rightarrow|n-m|\leq1.
$$
\end{defi}\label{def-comp}

\begin{figure}[ht]
\includegraphics*[width=11cm,height=4.5cm]{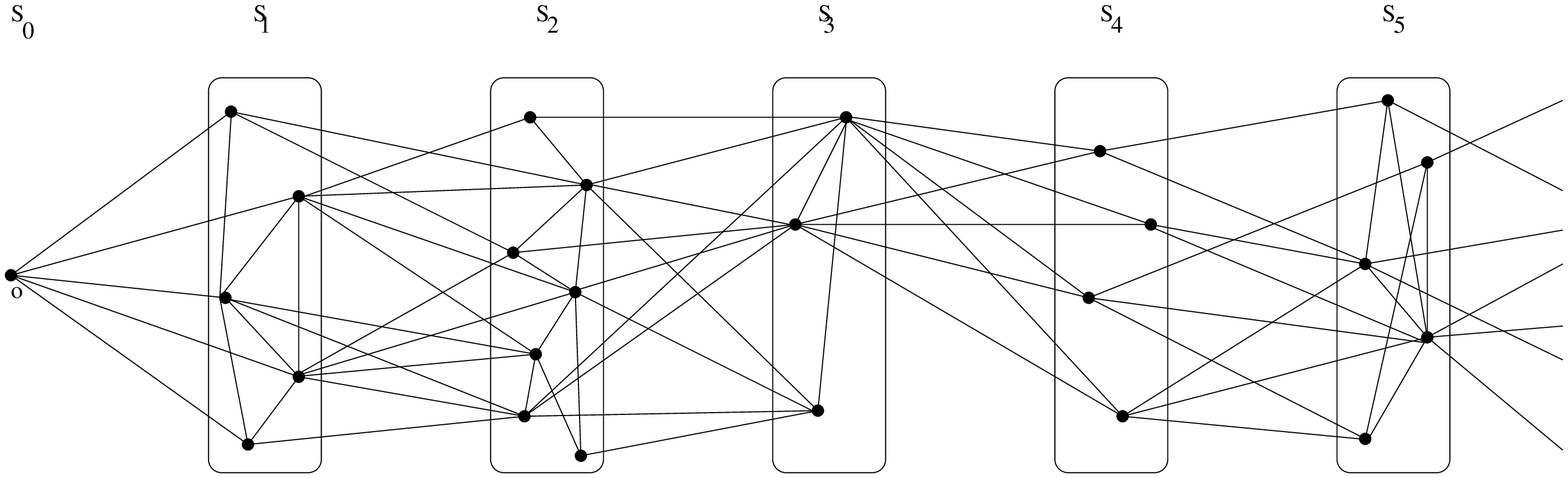}
\caption{\label{Paving} A triangulation with 1-dimensional decomposition}
\end{figure}

Given such a \emph{1-dimensional decomposition}, we write $B_n:=\cup_{i=0}^{n}\Sc_i.$ We set,
$$ (*)\left\{\begin{array}{ll}
       \mbox{deg}_{\Sc_n}^{\pm}(x):=\dfrac{1}{c(x)}
       \displaystyle\sum_{y\in\Vc(x)\cap\Sc_{n\pm1}}r(x,y) & \hbox{for all  }~ x\in\Sc_n, \\
       \mbox{deg}_{\Sc_n}^{0}(x):=\dfrac{1}{c(x)}\displaystyle\sum_{y\in\Vc(x)\cap\Sc_n}r(x,y)
       & \hbox{for all  }~ x\in\Sc_n, \\
       \mbox{deg}_{\Sc_n\times\Sc_{n+1}}(e):=\dfrac{1}{r(e)}
       \displaystyle\sum_{x\in\Fc_e\cap(\Sc_n\cup\Sc_{n+1})}s(e,x)
       & \hbox{for all  }~ e\in\Sc_n\times\Sc_{n+1}, \\
       \mbox{deg}_{\Sc_n^2}^0(e):=\dfrac{1}{r(e)}\displaystyle\sum_{x\in\Fc_e\cap\Sc_n}s(e,x)
       & \hbox{for all  }~ e\in\Sc_{n}^2, \\
       \mbox{deg}_{\Sc_n^2}^{\pm}(e):=\dfrac{1}{r(e)}\displaystyle\sum_{x\in\Fc_e\cap{\Sc_{n\pm1}}}s(e,x)
       & \hbox{for all  }~ e\in\Sc_{n}^2.\\

                   \end{array}
            \right.
   $$

We denote
$$\eta_n^{\pm}:=\displaystyle\sup_{x\in\Sc_n}\mbox{deg}_{\Sc_n}^{\pm}(x),~\beta_n:=\sup_{e\in\Sc_n\times\Sc_{n+1}}\mbox{deg}_{\Sc_n\times\Sc_{n+1}}(e),
~\gamma_n^{\pm}:=\sup_{e\in\Sc_n^2}\mbox{deg}_{\Sc_n^2}^{\pm}(e).
$$

\begin{thm}\label{thm.}
Let $\Tc=(\Kc,\Fc)$ be a triangulation and $(\Sc_n)_{n\in\N}$
a 1-dimensional decomposition of the graph $\Kc.$ Assume that
$$\displaystyle\sum_{n\in\N}\dfrac{1}{\sqrt{\xi(n,n+1)}}=\infty,
$$
with $\xi(n,n+1)=\eta_n^{+}+\eta_{n+1}^{-}+\beta_n+\gamma_n^{+}+\gamma_{n+1}^{-}.$
Then $\Tc$ is $\chi$-complete and in particular, $L$ is essentially self-adjoint
on $\Cc_{c}(\Vc)\oplus\Cc_{c}(\Ec)\oplus\Cc_{c}(\Fc).$
\end{thm}
\begin{dem}

We set
$$ \chi_n(x)=\left\{\begin{array}{ll}
                        1 & \hbox{if  }~ d_{comb}(o,x)\leqslant n, \\
                        \max\left(0,1-\displaystyle\sum_{k=n}^{d_{comb}(o,x)-1}
                        \dfrac{1}{\sqrt{\xi(k,k+1)}}\right)
                        & \hbox{if  } \, d_{comb}(o,x)> n.

                 \end{array}
            \right.
   $$

Since the series diverges, $\chi_n$ is with finite support.
Note that $\chi_n$ is constant on $\Sc_n.$ If $x\in\Sc_m$ with $m>n,$ we have
\begin{equation*}
\dfrac{1}{c(x)}\displaystyle\sum_{y\in\Vc(x)\cap\Sc_{m+1}}r(x,y)|\chi_n(x)-\chi_n(y)|^2
\leq\dfrac{\mbox{deg}_{\Sc_m}^{+}(x)}{\xi(m,m+1)}\leq1.
\end{equation*}
\begin{equation*}
\dfrac{1}{c(x)}\displaystyle\sum_{y\in\Vc(x)\cap\Sc_m}r(x,y)|\chi_n(x)-\chi_n(y)|^2=0.
\end{equation*}
\begin{equation*}
\dfrac{1}{c(x)}\displaystyle\sum_{y\in\Vc(x)\cap\Sc_{m-1}}r(x,y)|\chi_n(x)-\chi_n(y)|^2
\leq\dfrac{\mbox{deg}_{\Sc_m}^{-}(x)}{\xi(m-1,m)}\leq1.
\end{equation*}

On the other hand,
\begin{enumerate}
\item[-]If $e\in\Sc_m\times\Sc_{m+1},$ we get
\begin{equation*}
\dfrac{1}{r(e)}\displaystyle\sum_{x\in\Fc_e\cap(\Sc_m\cup\Sc_{m+1})}s(e,x)
|(\chi_n(x)-\chi_n(e^-))+(\chi_n(x)-\chi_n(e^+))|^2\leq
\dfrac{\mbox{deg}_{\Sc_m\times\Sc_{m+1}}^{+}(x,y)}{\xi(m,m+1)}\leq1.
\end{equation*}
\item[-]If $e\in\Sc^2_m,$ we get
\begin{equation*}
\dfrac{1}{r(e)}\displaystyle\sum_{x\in\Fc_e\cap\Sc_m}s(e,x)
|(\chi_n(x)-\chi_n(e^-))+(\chi_n(x)-\chi_n(e^+))|^2=0
\end{equation*}
\begin{equation*}
\dfrac{1}{r(e)}\displaystyle\sum_{x\in\Fc_e\cap\Sc_{m+1}}s(e,x)
|(\chi_n(x)-\chi_n(e^-))+(\chi_n(x)-\chi_n(e^+))|^2\leq
\dfrac{\mbox{deg}_{\Sc_m^2}^+(e)}{\xi(m,m+1)}\leq1.
\end{equation*}
\begin{equation*}
\dfrac{1}{r(e)}\displaystyle\sum_{x\in\Fc_e\cap\Sc_{m-1}}s(e,x)
|(\chi_n(x)-\chi_n(e^-))+(\chi_n(x)-\chi_n(e^+))|^2\leq
\dfrac{\mbox{deg}_{\Sc_m^2}^-(e)}{\xi(m-1,m)}\leq1.
\end{equation*}
\end{enumerate}
Then $\Tc$ is $\chi$-complete and in particular, $L$ is essentially self-adjoint by Corollary \ref{co5}.
\end{dem}
\subsection{A triangular tree}
Let $\Tc$ be a triangular tree, endowed with an origin. Due to its structure, one can take
$$ (**)\left\{\begin{array}{ll}
       \mbox{deg}_{\Sc_n}^{-}(x):=\dfrac{1}{c(x)}r(x,\overleftarrow{x})
       & \hbox{for all  }~ x\in\Sc_n, \\
       \mbox{deg}_{\Sc_n\times\Sc_{n+1}}(e)
       :=\dfrac{1}{r(e)}\displaystyle\sum_{x\in\Fc_e\cap\Sc_{n+1}}s(e,x)
       & \hbox{for all  }~ e\in\Sc_n\times\Sc_{n+1}, \\
       \mbox{deg}^-_{\Sc_n^2}(e):=\dfrac{1}{r(e)}s(e,\overleftarrow{e})
       & \hbox{for all  }~ e\in\Sc_{n}^2, \\

                   \end{array}
            \right.
   $$

where $\overleftarrow{e}$ is a unique vertex in $\Sc_{n-1}\cap\Fc_e.$
\begin{propr}
Let $\Tc$ be a triangular tree with its origin $o.$ Assume that
$$\displaystyle\sum_{n\in\N}\dfrac{1}{\sqrt{\xi(n,n+1)}}=\infty,
$$
with $\xi(n,n+1)=\eta^+_n+\eta^-_{n+1}+\beta_n+\gamma^-_{n+1}.$
Then $\Tc$ is $\chi$-complete and in particular, $L$ is essentially self-adjoint
on $\Cc_{c}(\Vc)\oplus\Cc_{c}(\Ec)\oplus\Cc_{c}(\Fc).$
\end{propr}
\begin{dem}
Use the same method as Theorem \ref{thm.} with $(**).$
\end{dem}

\subsection{Essential self-adjointness on the simple case}
In \cite{W} and \cite{D}, they prove that $L_0$ is essentially self-adjoint on $\Cc_{c}(\Vc)$
when the graph is simple. But the self-adjointness property does not always hold with other
operators in the simple case. We recall the operator $L_1^-$ is not necessarily essentially
self-adjoint on simple tree, see \cite{BGJ}. Moreover, we refer to \cite{GS} for the \emph{adjacency matrix}
$\Ac_{\Kc}=deg-L_0$ where $deg$ denotes the operator of multiplication with the functions which shows that
the deficiency indices of $\Ac_{\Kc}$ are infinite. In this framework, it is important to notice that $L_1$
and $L_2$ are not necessarily essentially self-adjoint on a simple triangulation.
\begin{propr}\label{pr}
Let $\Tc$ be a simple triangular tree. Assume that
\begin{equation}\label{eq*}
\emph{off}(n)=\#\left(\Vc(x)\cap\Sc_{n+1}\right),~x\in\Sc_n
\end{equation}
\begin{equation*}
n\mapsto\dfrac{\emph{off}^2(n)}{\emph{off}(n+1)}\in l^1(\N).
\end{equation*}

Then, $L_1$ is not essentially self-adjoint on $\Cc_{c}(\Ec)$
and the deficiency indices are infinite.
\end{propr}
\begin{dem}

We construct $\varphi\in l^2(\Ec)\backslash\{0\},$ such that
$\varphi\in Ker(L^*_1+i)$ and such that $\varphi$ is equal
to constant $C_n$ on $\Sc_n\times\Sc_{n+1}.$
It takes the value $0$ on $\Sc_n^2.$ Given the fact that $(x,y)\in\Sc_n^2,$ we get
$$C_n\left(\#\left(\Vc(x)\cap\Sc_{n+1}\right)-\#\left(\Vc(y)\cap\Sc_{n+1}\right)\right)=0.
$$

It holds because of the condition (\ref{eq*}). Now, we set $(x,y)\in\Sc_n\times\Sc_{n+1}$ and
with the condition (\ref{eq*}), we have
$$(\mbox{off}(n)+1+i)C_n-\mbox{off}(n+1)C_{n+1}-C_{n-1}=0.
$$

We can then apply Theorem 5.1 in \cite{BGJ} to obtain the conclusion. 
\end{dem}

We will see now also that $L_2$ is not necessarily essentially self-adjoint on simple
triangulation.
\begin{propr}\label{prr}
Let $\Tc$ be a simple triangulation with 1-dimensional decomposition
as shown in Figure $4.$ Assume that
\begin{equation}\label{tt}
n\mapsto\dfrac{\#\Sc_{2n}}{\#\Sc_{2(n+1)}}\in l^1(\N).
\end{equation}

Then, $L_2$ is not essentially self-adjoint on $\Cc_{c}(\Fc).$
\end{propr}
\begin{dem}
\begin{figure}[ht]
\includegraphics*[width=11cm,height=4.3cm]{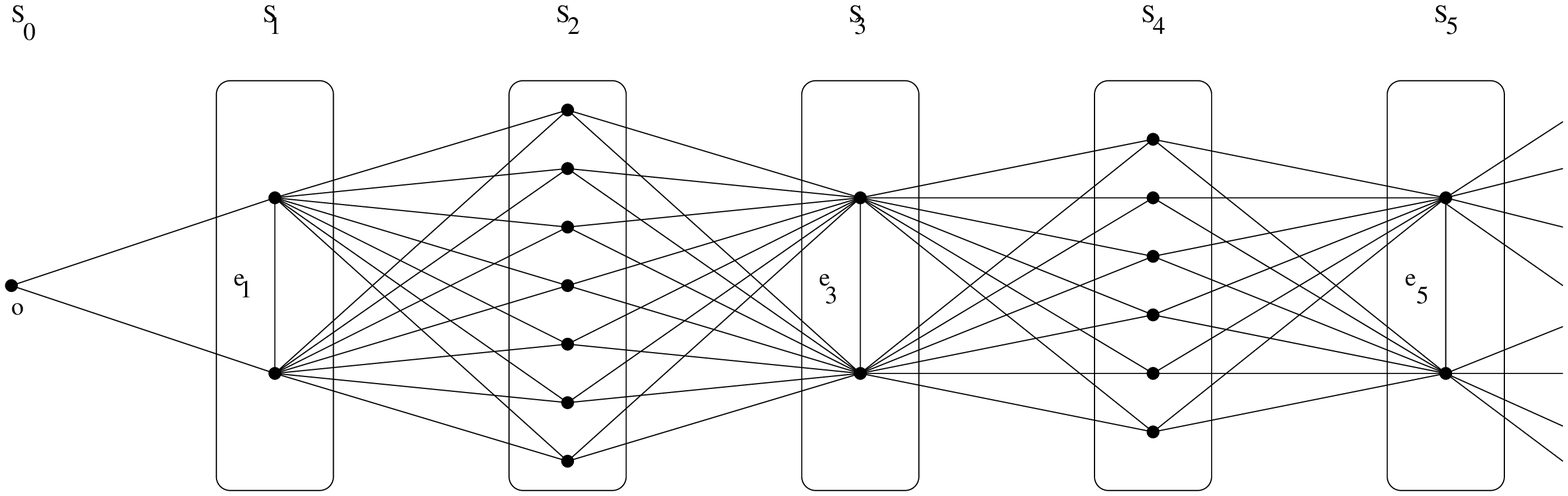}
\caption{\label{Paving} A particular triangulation with 1-dimensional decomposition}
\end{figure}

We consider the faces in Figure $4$ as follows:
$$\varpi\in\Fc\Leftrightarrow\mbox{ there exists }n\in\N\cup\{0\}\mbox{ such that } \varpi=(e_{2n+1},x)\in\left(\Sc_{2n+1}^2\times\Sc_{2n}\right)
\cup\left(\Sc_{2n+1}^2\times\Sc_{2n+2}\right).
$$

Set $\phi\in l^2(\Fc)\backslash\{0\}$ such that $\phi\in Ker(L^*_2+i).$
For $n\in\N,$ it is given by
$$ \phi(e_{2n+1},x):=\left\{\begin{array}{ll}
       C_{2n+2} & \hbox{for all  }~ x\in\Sc_{2n+2}.\\
       C_{2n} & \hbox{for all  }~ x\in\Sc_{2n}.\\

                   \end{array}
            \right.
   $$

Let $x\in\Sc_{2n+2},$ we have
\begin{equation*}
\begin{split}
(L^*_2+i)(\phi)(e_{2n+1},x)&=\displaystyle\sum_{u\in\Fc_{e_{2n+1}}}\phi(e_{2n+1},u)
+\displaystyle\sum_{u\in\Fc_{(e^+_{2n+1},x)}}\phi(e^+_{2n+1},x,u)
+\displaystyle\sum_{u\in\Fc_{(x,e^-_{2n+1})}}\phi(x,e^-_{2n+1},u)\\
&+i\phi(e_{2n+1},x)=0.
\end{split}
\end{equation*}

Hence, we get
\begin{equation}\label{et}
\left(\#\Sc_{2n+2}+2+i\right)C_{2n+2}+\left(\#\Sc_{2n}\right)C_{2n}=0.
\end{equation}

By the equation (\ref{et}), we get
\begin{equation*}
\begin{split}
\|\phi_{\scriptscriptstyle{\vert \Sc_{2n+1}^2\times\Sc_{2n+2}}}\|_{l^2(\Fc)}^2
&=\dfrac{1}{6}\displaystyle\sum_{[x,y,z]\in\Sc_{2n+1}^2\times\Sc_{2n+2}}|\phi(x,y,z)|^2\\
&=\dfrac{1}{6}\left(\#C_{2n+2}\right)^2\left(\#\Sc_{2n+2}\right)\\
&=\dfrac{1}{6}\dfrac{\left(\#C_{2n}\right)^2
\left(\#\Sc_{2n}\right)^2}{|\#\Sc_{2n+2}+2+i|^2}\left(\#\Sc_{2n+2}\right)\\
&=\dfrac{\left(\#\Sc_{2n}\right)\left(\#\Sc_{2n+2}\right)}{|\#\Sc_{2n+2}+2+i|^2}
\|\phi_{\scriptscriptstyle{\vert \Sc_{2n-1}^2\times\Sc_{2n}}}\|_{l^2(\Fc)}^2
\end{split}
\end{equation*}

Since $\displaystyle\lim_{n\rightarrow\infty}\dfrac{\#\Sc_{2n}}
{\#\Sc_{2(n+1)}}=0,$ we get by induction
$$C:=\displaystyle\sup_{n\in\N^*}\|\phi_{\scriptscriptstyle{\vert \Sc_{2n-1}^2\times\Sc_{2n}}}
\|_{l^2(\Fc)}^2<\infty.
$$

Thus
$$\|\phi_{\scriptscriptstyle{\vert \Sc_{2n+1}^2\times\Sc_{2n+2}}}\|_{l^2(\Fc)}^2
\leq C\dfrac{\left(\#\Sc_{2n}\right)\left(\#\Sc_{2n+2}\right)}{|\#\Sc_{2n+2}+2+i|^2}.
$$

From (\ref{tt}), we conclude that $\phi\in l^2(\Fc).$ By mimicking the proof of Theorem X.36 of
\cite{RSv2} one shows that $\mbox{dimKer}(L_2^*+i)\geq1$ and thus $L_2$ is not essentially
self-adjoint on $\Cc_{c}(\Fc).$
\end{dem}
\begin{rem}
By one of Propositions \ref{prr} and \ref{pr}, we conclude that $L$ is not necessarily essentially
self-adjoint on a simple triangulation.
\end{rem}

\textbf{\textit{Acknowledgments}:} I would like to sincerely thank my PhD advisors,
\emph{Colette Ann\'e} and \emph{Nabila Torki-Hamza} for helpful discussions.
I am very thankful to them for all the encouragement, advice and inspiration.
I take this chance to thank Matthias Keller for the fruitful discussions during my visit to Potsdam.
I would also like to thank the Laboratory of Mathematics Jean Leray and the research unity
(UR/13 ES 47) for their continuous support. This work was financially supported by the "PHC Utique"
program of the French Ministry of Foreign Affairs and Ministry of higher education and research
and the Tunisian Ministry of higher education and scientific research in the CMCU project number 13G1501
"Graphes, G\'{e}om\'{e}trie et Th\'{e}orie Spectrale". Finally, I would like to thank the referee for
the careful reading of my paper and the valuable comments and suggestions.

\end{document}